\documentclass[11pt]{amsart}
\usepackage{amsmath,amsthm,amsfonts,amssymb}
\usepackage{amsmath, amsthm, amssymb, amscd, amsfonts}

\usepackage[all]{xy}
\usepackage{amssymb, amsthm, amsmath}
\usepackage[english]{babel}
\usepackage{amsfonts}
\usepackage{color}
\usepackage{leftidx}
\usepackage{enumerate}
\usepackage[dvipsnames]{xcolor}

\usepackage{mathtools}
\usepackage{mathtools}
\usepackage[english]{babel}
\usepackage{tikz}
\usetikzlibrary{matrix}
\newcommand{\re}{\res}
\usepackage{leftidx}
\usepackage{enumerate}

\newcommand{\bll}{\blue}

\DeclareMathOperator{\id}{id}
\DeclareMathOperator{\Hom}{Hom}

\newcommand{\ra}{\rightarrow}

\newcommand{\Z}{\mathbb Z}
\newcommand{\ot}{\otimes}
\newcommand{\mtc}{\mathcal}

\newcommand{\al}{\alpha}
\newcommand{\eps}{\epsilon}
\newcommand{\bn}{\begin}

\newcommand{\cn}{\mathcal{N}}

\newcommand{\sub}{\subsection}

\newcommand{\D}{\Delta}

\newcommand{\ct}{\mathcal M}

\numberwithin{equation}{section}
\newtheorem{lemma}[equation]{Lemma}

\newtheorem{defn}[equation]{Definition}
\newtheorem{cor}[equation]{Corollary}
\newtheorem{rem}[equation]{Remark}

\newcommand{\bl}{\begin{lemma}
  }
  \newcommand{\M}{\mathcal{M}}
\newcommand{\un}{\mathbf{1}}
\newcommand{\nc}{\newcommand}
\nc{\el}{\end{lemma}}

\newcommand{\ch}{\chi}
\newcommand{\mtr}{\mathrm}
\nc{\bwt}{\bowtie}
\newcommand{\ncm}{\newcommand}\newcommand{\gm}{\gamma}
\numberwithin{equation}{section}
\newcommand{\et}{\end{thm}}\newcommand{\bt}{\bn{thm}}
\newcommand{\ep}{\end{prop}}\newcommand{\bp}{\bn{prop}}
\newcommand{\beqarn}{\begin{eqnarray*}}
\newcommand{\eeqarn}{\end{eqnarray*}}
\newcommand{\beqn}{\bn{equation*}}
\newcommand{\eeqn}{\end{equation*}}
\newcommand{\bpf}{\bn{proof}}
\newcommand{\epf}{\end{proof}}
\ncm{\cX}{\mtc{X}}
\ncm{\wt}{\widetilde}
\ncm{\sg}{\sigma}\ncm{\Rep}{\mathrm{Rep}}
\newcommand{\Res}{\mathrm{Res}}
\ncm{\X}{\mathcal{X}}
\ncm{\cA}{\mathcal{A}}
\newcommand{\lb}{\label}

\numberwithin{equation}{section}

\numberwithin{equation}{section}

\ncm{\np}{\newpage}
\ncm{\ebl}{\end{thebibliography}}
\ncm{\bbl}{\begin{thebibliography}}
\ncm{\chd}{_{ _{\ch}}}
\ncm{\ald}{_{ _{\al}}}
\ncm{\cP}{\mathcal{P}}
\ncm{\ei}{e_i}
\ncm{\eij}{e_{i,\;j}}
\ncm{\bne}{\begin{enumerate}}
\ncm{\ene}{\end{enumerate}}\ncm{\bdef}{\begin{defn}}
\ncm{\edf}{\end{defn}}
\ncm{\stab}{\mtr{Stab}}
\ncm{\bc}{\begin{cor}}

\ncm{\ec}{\end{cor}}
\ncm{\er}{\end{rem}}
\ncm{\br}{\begin{rem}}

\ncm{\bd}{
\begin{document}}
\ncm{\ed}{\end{document}}
\ncm{\beq}{\begin{equation}}

\ncm{\eeq}{\end{equation}}
\ncm{\cm}{\mathcal{M}}
\ncm{\rep}{\mtr{Rep}}
\ncm{\btw}{\bowtie}
\ncm{\cd}{\mtc{D}}
\ncm{\cop}{\mtr{cop}}

\ncm{\bea}{\begin{eqnarray}}
\ncm{\eea}{\end{eqnarray}}
\ncm{\beanon}{\begin{eqnarray*}}
\ncm{\eeanon}{\end{eqnarray*}}\ncm{\ek}{\eps|_K}\ncm{\diez}{\#}

\ncm{\cC}{\mtc{C}}
\ncm{\cc}{\mtc{C}}
\ncm{\HKer}{\mtr{HKer}}
\ncm{\LKER}{\mtr{LKER}}
\ncm{\aad}{\mtr{ad}}
\ncm{\Dr}{\mtr{D}}
\ncm{\cD}{\mathcal{D}}
\ncm{\G}{\mathcal{G}}
\ncm{\Dc}{\mtc{D}}
\ncm{\E}{\mtc{E}}
\ncm{\fp}{\mtr{FP}}
\ncm{\Vc}{\mtr{Vec}}
\ncm{\cK}{\mtc{K}}
\ncm{\cM}{\mtc{M}}
\ncm{\cE}{\mtc{E}}
\ncm{\cS}{\mtc{S}}
\ncm{\cs}{\cS}
\ncm{\End}{\mtr{End}}
\ncm{\hsa}{Hopf subalgebra of }
\ncm{\ses}{semisimple}
\ncm{\x}{$}
\newcommand{\mdn}{\medskip\noindent}
\ncm{\mi}{\mtr{I}}
\ncm{\cZ}{\mtc{Z}}\ncm{\xra}{\xrightarrow}
\ncm{\cb}{\mtc{B}}\ncm{\ca}{\mtc{A}}
\ncm{\irr}{\Irr}\ncm{\Irr}{\mathrm{Irr}}
\ncm{\co}{\mtc{O}}
\ncm{\cg}{{\mtr{K}_0}}\ncm{\ci}{\mtc{I}}
\ncm{\blue}{\textcolor[rgb]{.00, .00, 1.00}}
\ncm{\bb}{\blue}
\ncm{\red}{\textcolor[rgb]{1.00, .00, .00}}
\ncm{\md}{\medbreak}
\ncm{\green}{\textcolor[rgb]{.00, 1.00, .00}}
\ncm{\Gm}{\Gamma}\ncm{\ind}{\mtr{Ind}}\ncm{\res}{\mtr{Res}}
\numberwithin{equation}{section}
\ncm{\bq}{\beq}\ncm{\mto}{\mapsto}\ncm{\opl}{\oplus}
\ncm{\eq}{\eeq}\newcommand{\R}{{\mathcal R}}
\newcommand{\C}{{\mathcal C}}
\ncm{\Ind}{\mtr{Ind}}\ncm{\cz}{\mtc{Z}}\ncm{\ce}{\mtc{E}}\ncm{\ro}{T}
\ncm{\inv}{^{{-1}}}
\title[Fusion categories]
{Categorical Green functors arising from group actions on categories}
\author{Sebastian  Burciu}
\address{Inst.\ of Math.\ ``Simion Stoilow" of the Romanian Academy, Research unit 5,
P.O. Box 1-764, RO-014700, Bucharest, Romania}\email{sebastian.burciu@imar.ro} 
\thanks{ The research was done in part during a stay at the Erwin Schr\"odinger Institute, Vienna, in the frame of the Programme ÒModern Trends in Topological Quantum
Field TheoryÓ in February 2014. The author thanks the ESI and the organizers of the Programme for the support and kind hospitality.
\medbreak
This work was supported by a grant of the Romanian National Authority for Scientific Research, CNCS-UEFISCDI, project number PN-II-RU-TE-2012-3-0168.}

\subjclass[2000]{Primary 16T20, 18D10, 19A22}
\keywords{Group actions on categories; Tensor categories;K-theory; Grothendieck rings;}
\bd
\maketitle

\begin{abstract}
In this paper we introduce the notion of a categorical  Mackey functor.  This categorical notion allows us to obtain new Mackey functors by passing to Quillen's $K$-theory of the corresponding  abelian categories. In the case of an action by monoidal autoequivalences on a monoidal category the Mackey functor obtained at the level of Grothendieck rings has in fact a Green functor structure. 
\end{abstract}
\section{Introduction and Main Results}
A Mackey functor or (a $G$-functor) is a family $\{a(K)\}_{K\leq G}$ of abelian groups equipped
with three types of maps: induction, conjugation, and restriction, satisfying some certain compatibility axioms, see for example \cite{gr71}. Typical examples, include among others, the cohomology groups $\{H^{n}(K, M)\}_{K\leq G}$ and the character rings $\{R(K)\}_{K \leq G}$.

In \cite{8, 35, 50} it is shown that for any group $G$ of automorphisms of a number field $k$, the class group of the ring of integers of the fixed field $\{k^{H}\}_{H \leq G}$ is a $G$-functor.
These results were extended in \cite{kevin} by showing that $\{K_{i}(S^{H })\}_{H\leq G}$ is a $G$-functor, whenever $R\subseteq S$ is a Galois extension of commutative rings with Galois group $G$.

The main goal of this paper is to construct a categorical version of a Mackey functor.  The main source of examples for such categorical functors is given by the group actions on categories. It is shown that group actions on abelian category give rise to Mackey functors while monoidal group actions on monoidal categories give rise to categorical Green functors. By passing to the K-theory a categorical Mackey functor give rise to a classical Mackey functor. In this way we give new examples of Mackey and Green functors generalizing the examples given \cite {kevin, jlms}.

Let $G$ be a finite group acting on the abelian category $\cc$. For any subgroup $H$ of $G$ the left adjoint functor of the forgetful functor $\res^{G}_{H}:\cc^{G}\ra  \cc^{H}$ was recently described in \cite{buna}.  This functor is denoted by $\Ind^{G}_{H}:\cc^{H}\ra  \cc^{G}$ and can be regarded as a generalization of the induction functor from $\rep(H)$ to $\rep(G)$.

Using this notion of induction and restriction we introduce the concept of categorical Mackey and Green functors that category the classical concepts. Note that the induction and restriction functors were also considered in various other contexts as categorized constructions of representation theory, see for example \cite{rq}.

\medbreak
Our first main result is the following: 
 \bt\lb{cmain}
Let $G$ be a finite group acting on the abelian category $\cc$ by $T:G \ra \underline{\mathrm{Aut}}(\cc)$. 
\bne
\item Then  the functor $H \mapsto \cc^{H}$ defines a categorical Mackey functor over $\mathrm{Vec}_{k}$. \item Moreover if $\cc$ is a $k$-linear monoidal category and the action of $G$ is by monoidal autoequivalences then the above functor is a categorical Green functor over $\cc^{G}$.  \ene
\et
\noindent
The proof of the above results uses a Mackey type decomposition for the above induced functor when restricted to various subgroups:
\bt\lb{macky}
Suppose that a finite group $G$ acts on the abelian category 
$\cc$ via $T:\underline{G} \to \underline{ \mtr{Aut}}(\cc)$. Let $K$ and $L$ be any two subgroups of a subgroup  $H \leq G$ and $M\in \cc^H$.

\noindent 1) Then
\beq
\res_K^H(\ind^H_L(M))\simeq\bigoplus_{x \in K\backslash H\slash L}\ind_{K \cap \;^xL}^K(\res^{\;^xL}_{\;^xL\cap K}(c_{L,x}(M)))
\eeq
where $\;^{x}L:=xLx^{{-1}}$ and the equivariant structure for $c_{L, x}(M)\in \cc^{\;^{x}L}$ is given as in Lemma \ref{conj}.

\noindent 2) If $\cc$ is a $k$-linear monoidal category  and the action of $G$ on $\cc$ is by monoidal autoequivalences then the above isomorphism is of $\cc^{G}$-module functors.
\et
As an application of Theorem \ref{kpass} we obtain the following corollary:
\bc\label{main3}
Let $G$ be a finite group acting on the abelian category $\cc$ by $T:G \ra \underline{\mathrm{Aut}}(\cc)$. 
\bne 
\item Then for all $i \geq 0$ the functor $H \mapsto K_{i}(\cc^{H}) $ defines a $G$-functor $M_{i}$ with the following structure maps:
\bne
\item Restriction
$
R^{H}_{K} :K_{i}(\cc^{H}) \ra K_{i}(\cc^{K})
$
is the map induced by the forgetful functor $\res^{H}_{K}:
\cc^{H}\ra \cc^{K}$, 
\item Induction $
I^{H}_{K} :K_{i}(\cc^{H}) \ra K_{i}(\cc^{K})
$
is the map induced by the induction functor $\ind^{H}_{K}:\cc^{K}\ra \cc^{H}$,
\item Conjugation $
c_{H, x} :K_{i}(\cc^{H})\ra K_{i}(\cc^{\;^{x}H})
$
is the map induced by the functor $T^{x}:\cc^{H}\ra \cc^{\;^{x}H}$.
\ene
\item If $\cc$ is a $k$-linear monoidal category  over $k$ and $G$ is a finite group acting  on $\cc$ by monoidal autoequivalences then $H \mapsto K_{0}(\cc^{H})$ defines a Green functor on $G$ over $k$.
\ene
\ec
\noindent
Shortly, this paper is organized as follows.  In Section \ref{prelim} we recall some basic results on abelian categories and group actions on them. The construction of the adjoint functor $\Ind^{G}_{H}$ mentioned above is also recalled in this section.

In Section \ref{catdef} we first recall the definition of the classical Mackey and Green functors.  Then we present the new categorical notions of Mackey and Green functors.

Section \ref{tim} is devoted to the proof of Theorem \ref{cmain} and Theorem \ref{macky}. \noindent In Theorem \ref{kpass} we show that by passing to Quillen's K-theory a categorical Mackey functor gives rise to a classical Mackey functors. 
(see Theorem \ref{kpass}). This proves Corollary \ref{main3}.
\section{Group actions on $k$-linear categories}\lb{prelim}
Fix a commutative ring $k$. Recall that a $k$-linear category is an abelian category in which the hom-sets are $k$-vector spaces and the compositions of morphisms are $k$-bilinear. A $k$-linear functor between $k$-linear categories is a functor which is linear on all hom-spaces. Recall that an adjunction between categories C and D is a pair of functors,
$
F: \mathcal{D} \rightarrow \mathcal{C}$   and   $G: \mathcal{C} \rightarrow \mathcal{D}$
and a family of bijections
$
\mathrm{hom}_{\mathcal{C}}(FY,X) \cong \mathrm{hom}_{\mathcal{D}}(Y,GX)$
which is natural in the variables $X$ and $Y$. The functor F is called a left adjoint functor, while G is called a right adjoint functor. The relationship ÒF is left adjoint to GÓ (or equivalently, ÒG is right adjoint to FÓ) is sometimes written $F\dashv G$.

\noindent
A {\it monoidal category} is a category $\cc $ equipped with
a bifunctor $\otimes \colon \cc\times\cc\to\cc$ called the {\it  monoidal product}, and an object $\un_{\cc}$ called the {\it unit object}.  The category $\cc $ has a natural isomorphism $\alpha$, called {\it associativity constraint}, given by $\alpha_{A,B,C}  (A\otimes B)\otimes C \xra{\simeq} A\otimes(B\otimes C)$ for all $A, B, C \in\cc$. There are also two natural isomorphisms $l_{A}$ and $r_{A}$, respectively called left and  respectively {\it right unitor}, with components $l_A \colon I\otimes A\cong A$ and $r_A \colon A\otimes I\cong A$. 
These natural transformations satisfy some coherence conditions expressed by the fact that pentagon and the  triangle diagram commute, see e.g. \cite{ENO}.

Recall that a {\it unitary monoidal functor} $F:\cc \ra \cd$ between two monoidal categories is a $k$-linear functor $F$ together with a natural transformation $F_{2}:F(- \; \ot \; -)\ra F(-)\ot F(-)$ and a unit isomorphism $F_{0}:F(\un_{\cc}) \ra \un_{\cd}$ satisfying the compatibility of the following  hexagon and unit axioms (see for example \cite{ENO}). {
\begin{center}
{\tiny
\begin{equation}\tag{Diagram H}
\begin{tikzpicture}
  \matrix (m) [matrix of math nodes, row sep=6.5em,column sep=6em,minimum width=5em]
  {  (F(A)\ot F(B))\ot F(C) &  F(A)\ot (F(B)\ot F(C)) \\
 F(A  \ot B) \ot F(C) & F(A) \ot F(B\ot C)\\
  F((A\ot B)\ot C) &  F(A\ot (B\ot C)).  \\
    };
  \path[-stealth]
    (m-1-1) edge node [above] {$\al_{F(A), F(B), F(C)}$} (m-1-2)
    (m-3-1) edge node [above] {$F(\al_{A, B, C})$} (m-3-2)
    (m-1-1) edge node [left] {$ F_{2}^{A, B}\ot 1$} (m-2-1)
    (m-1-2) edge node [right] {$1\ot F_{2}^{B, C}$} (m-2-2)
    (m-2-1) edge node [left] {$ F_{2}^{A\ot B, C}$} (m-3-1)
    (m-2-2) edge node [right] {$F_{2}^{A, B\ot C}$} (m-3-2)
   ;
\end{tikzpicture}
\end{equation}
}
\end{center}
}
{
\begin{center}
{\tiny
\begin{equation}\tag{Diagram LU}
\begin{tikzpicture}
  \matrix (m) [matrix of math nodes, row sep=3em,column sep=3em,minimum width=0.2em]
  {  F(\un_{\cc} \ot A)  & F(A) \\
  F(\un_{\cc} ) \ot F( A) & \un_{\cd}  \ot F( A)\\
    };
  \path[-stealth]
    (m-1-1) edge node [left] {$F_{2}^{\un_{\cc}, A}$} (m-2-1)
    (m-2-1) edge node [below] {$F_{0}\ot 1$} (m-2-2)
    (m-1-1) edge node [above] {$ F(l_{A})$} (m-1-2)
    (m-2-2) edge node [right] {$l_{F(A)}$} (m-1-2)
   ;
\end{tikzpicture}
\end{equation}
}
\end{center}
}
\begin{center}
{\tiny
\begin{equation}\tag{Diagram RU}
\begin{tikzpicture}
  \matrix (m) [matrix of math nodes, row sep=3em,column sep=3em,minimum width=0.2em]
  {  F(A \ot \un_{\cc} )  & F(A) \\
 F(A)\ot F(\un_{\cc}) & F( A) \ot \un_{\cd}  \\
    };
  \path[-stealth]
    (m-1-1) edge node [left] {$F_{2}^{A, \un_{\cc}}$} (m-2-1)
    (m-1-1) edge node [above] {$ F(r_{A})$} (m-1-2)
    (m-2-1) edge node [below] {$1 \ot F_{0}$} (m-2-2)
    (m-2-2) edge node [right] {$r_{F(A)}$} (m-1-2)
   ;
\end{tikzpicture}
\end{equation}
}
\end{center}
\noindent
In particular the naturally of $F_{2}$ with respect to the morphisms can be written as
\beq\label{natf}
(F(u)\ot F(v))F_{2}^{M, N}=F_{2}^{M', N'}F(u \ot v)
\eeq
for all morphisms $M\xra{u} M'$ and $N\xra{v} N'$ in $\cc$.
Composition of two monoidal functors $\cc \xra{G} \cd \xra{F}\ce$ is also a monoidal functor with 
\beq\label{comp}
(F\circ G)_{2}^{M, N}:=F_{2}^{G(M), G(N)}\circ F(G_{2}^{M, N})
\eeq
A {\it natural monoidal transformation }$\tau: F\ra G$ between two monoidal functors is a natural transformation satisfying the following compatibility condition:
\beq\label{tenant}
G^{M, N}_{2}\tau_{M\ot N}=(\tau_{M}\ot \tau_{N})F^{M, N}_{2}
\eeq
for any objects $M, N \in \cc$. 

\noindent
If $\C$ is a monoidal category, a {\it left module category} over $\C$ is a  category $\M$ endowed with an \emph{action} functor $\boxtimes : \C \times \M \to \M$ with  {\it module associativity constraints} $a_{A, B, M}: A \boxtimes (B \boxtimes M)\xra{\simeq} (A \otimes B) \boxtimes M$ and a unit constraint $u_{M}:\un \boxtimes M \xra{\simeq} M$ satisfying a pentagon and a triangle coherence axiom, see \cite{ENO}. A {\it module functor} $F:\cm \ra \cn$ between two $\cc$-module categories $\cm,\cn$ is a functor with an additional module structure such that
$
(F_{2})^{A, M}:F(A\boxtimes M)\simeq A \boxtimes F(M)
$
satisfying the following diagrams: 
\newcommand{\bxt}{\boxtimes}
{\tiny
\begin{center}
\begin{equation}\tag{Diagram UM}
\begin{tikzpicture}
  \matrix (m) [matrix of math nodes, row sep=4em,column sep=6em,minimum width=5em]
  {  F(\un \boxtimes M) & & F(M) \\
 & \un  \bxt F( M) & \\
    };
  \path[-stealth]
    (m-1-1) edge node [left] {$F_{2}^{\un, M}$} (m-2-2)
    (m-1-1) edge node [above] {$ F(u_{M})$} (m-1-3)
    (m-2-2) edge node [right] {$u_{F(M)}$} (m-1-3)
   ;
\end{tikzpicture}
\end{equation}
\end{center}
}
and {\tiny
\begin{center}
\begin{equation}\tag{Diagram HM}
\begin{tikzpicture}
  \matrix (m) [matrix of math nodes, row sep=4em,column sep=6em,minimum width=5em]
  {  F((A\ot B)\boxtimes M) &  & (A\ot B)\bxt F(M) \\
 F(A  \bxt (B \bxt M)) & & A \bxt (B\bxt F(M))\\
 & A\bxt F(B\ot M) &    \\
    };
  \path[-stealth]
    (m-1-1) edge node [above] {$F_{2}^{A\ot B, M}$} (m-1-3)
    (m-3-2) edge node [below] {$1\ot F_{2}^{B,M}$} (m-2-3)
    (m-1-1) edge node [left] {$ F(\al_{A, B, F(M)})$} (m-2-1)
    (m-1-3) edge node [right] {$\al_{A, B, F(M)}$} (m-2-3)
    (m-2-1) edge node [left] {$ F_{2}^{A, B\bxt M}$} (m-3-2)
   ;
\end{tikzpicture}
\end{equation}
\end{center}
}
\noindent Let $\cc$ be a monoidal category and $\cm$, $\cn$ be two left $\cc$-module categories. Let also $F, G:\cm \ra \cn$ be two $\cc$-module functors. A morphism between $F
$ and $G$ is a natural transformation $T:F\ra G$ such that the following diagram commutes
\begin{center}
{\tiny
\begin{tikzpicture}
  \matrix (m) [matrix of math nodes, row sep=6em,column sep=6em,minimum width=5em]
  { F_{1}(P \otimes M) & F_{2}(P \ot M) \\
   P \otimes F_{1}( M)  &   P \otimes F_{2}( M)\\
    };
  \path[-stealth]
  (m-1-1) edge node [above] {$T_{P\ot M}$} (m-1-2)
    (m-2-1) edge node [above] {$1_{P}\ot T_{M}$} (m-2-2)
         
               (m-1-1) edge node [left] {$F_{2}^{P,M}$} (m-2-1)
    (m-1-2) edge node [right] {$G_{2}^{P,M}$} (m-2-2)
   ;
\end{tikzpicture}
	}
   \end{center}
   for any $P \in \cc$ and any $M \in \cm$.
\noindent
 Let $R, S:\cc\ra\cd$  be two functors and $N:R\ra S$ be a natural transformation between $R$ and $S$. For any functor $F:\cd \ra \ce$ one can define the natural transformation $_{F}N:FR\ra FS$ by $(_{F}N)_{X}:=F(N_{X}):FR(X)\ra FS(X)$.

\noindent
Moreover for any functor $G:\ce \ra \cc$ one can also define $N_{G}:RG \ra SG$ as natural transformation  by $(N_{G})_{X}:=N_{G(X)}$. 
\subsection{Group actions on abelian categories}
Let $\C$ be an abelian category. Denote by $ \underline{ \mtr{Aut}}(\cc)$ the category whose objects are {\it exact} autoequivalences of $\cc$ and morphisms are natural transformations between them. Then $ \underline{ \mtr{Aut}}(\cc)$ is a monoidal category where the tensor product is defined as the composition of autoequivalences. For a finite group $G$ let $\mathrm{Cat}(G)$ denote the monoidal category whose objects
are elements of $G$, the only morphisms are the identities, and the tensor product is given by multiplication in $G$. An action of a finite group $G$ on $\cc$ consists of a unitary monoidal functor $T:\mathrm{Cat}(G) \ra \underline{ \mtr{Aut}}(\cc)$. Thus, for every $g \in G$, we have a functor $T^g: \C \to \C$ and a collection of natural isomorphisms 
\beqn T^{g,h}_2 : T^g T^h \to T^{gh}, \quad g, h \in G,\eeqn which give the monoidal structure of $T$. The monoidal unit of $T$ is denoted by $T_{0}:\id_{{\cc}}\ra T^{1}$ where $1 \in G$ is the unit of the group $G$.  By the definition of the monoidal functor, the monoidal structure $T_{2}$ satisfies the following conditions:
\beq\label{ro-2}  (T^{gh, l}_2)_M \, (T^{g,h}_2)_{T^l(M)} =
(T^{g, hl}_2)_M \, T^g((T^{h,l}_2)_M), \eeq\beq \label{gunit}  (T^{g,
1}_2)_M T^{g}({T_0}_{ _{M}}) = (T^{1, g}_2)_M
(T_0)_{T^{g}(M)}=\id_{T^{g}(M)},
\eeq
for all objects $M \in \C$, and for all $g, h, l \in G$. See \cite[Subsection 4.1]{DGNO}. 
Note that by the naturality of $T^{g, h}_2$, $g, h \in G$,  can be written as
\begin{equation}\label{nat-ro}T^{gh}(f) \, (T_2^{g, h})_{N} = (T_2^{g, h})_M \, T^gT^h(f),
\end{equation}
 for every morphism $f: M \to N$ in $\C$. {\it 
We shall assume in what follows that $T^1 = \id_\C$ and
$T_0$, $T^{g, 1}_2$, $T^{1, g}_2$ are also identities.}  We say that $G$ acts {\it $k$-linearly} on the $k$-linear category $\cc$ if $T^{g}$ is a $k$-linear autoequivalence for any $g \in G$.
\bn{example}\label{automact}
Suppose that $G$ acts by ring automorphisms on a $k$-algebra $S$. Then $G$ acts on $S$-mod via  the following action: $T^{g}(M)=M$ as abelian groups and the $S$-action on $T^{g}(M)$ is given by $s.\;^{g}m:=(g^{-1}.s)m$. In this case one can take $(T^{g, h}_{2})_{M}=\id_{M}$ for all $g,h \in G$.
\end{example}

\subsection{On the equivariantized category}
\medbreak Suppose that $G$ acts on the abelian category $\cc$. Let $\C^G$ denote the corresponding
\emph{equivariantized} category. Recall that $\C^G$ is an abelian category whose objects are $G$-equivariant
objects of $\C$. They consist of pairs $(M, \mu)$, where $M$ is an object
of $\C$ and $\mu = (\mu_{M}^g)_{g \in G}$ is a collection of isomorphisms $\mu_{M}^g:T^g(M) \to M$ in $\cc$ satisfying the following:
\begin{equation}\label{deltau} \mu_{M}^g T^g(\mu_{M}^h) = \mu_{M}^{gh} (T^{g,
h}_2)_M, \quad \forall g, h \in G, \quad \quad
\mu_{M}^{1}{T_0}_M=\id_M.\end{equation} \medbreak We say that an object $M$ of $\C$ is
\emph{$G$-equivariant} if there exists such a collection $\mu =
(\mu^g)_{g \in G}$ so that $(M, \mu) \in \C^G$. Note that the equivariant structure $\mu$ is not necessarily unique.

\noindent
A morphism $f : (M, \mu_{M})\ra (N, \mu_{N})$ in $\cc^{G}$ is a morphism in $f:M\ra N$ such that 
\beq\label{morf}
\mu^{g}_{N}T^{g}(f)=f\mu^{g}_{M}.
\eeq
\bn{example}\label{automequiv}
It is easy to verify that in the case of the previous example one has that $(S$-mod)$^{G}\simeq S\#kG$-mod, the category of $S\#kG$-modules.
\end{example}
\subsection{Induction functors as left adjoints of restriction functors}\label{indc} Suppose that a finite group $G$ acts on the abelian category $\cc$ and let $H\leq G$ be a subgroup. Let $\R$ be a set of representative elements for the left  cosets $\{Hx\;x \in G\}$
of $H$ in $G$. Thus one can write $G$ as a disjoint union $G = \cup_{t \in \R}tH$. 
Set, for all $(V, \mu) \in \cc^{H}$, 
\bq\lb{indx}
\ind_H^G(V, \mu): =(\oplus_{t \in \R}T^t(V), \nu)\in \cc^{G}
\eq 
where for all $g \in G$ the equivariant structure of $\nu^{g}_{\ind_H^G(V, \mu)}: \bigoplus_{t \in \R}T^gT^t(V) \to \bigoplus_{t \in \R}T^t(V)$ is defined componentwise by the formula 
\begin{equation}\label{muind}
\nu^{g, t} = T^s(\nu^h)
(T_2^{s, h})^{-1} T_2^{g, t}: T^g T^t(V) \to
T^s(V).\end{equation} 
Here the elements $h \in H$ and $s \in \R$ are uniquely determined by the relation $gt = sh$.
\bp \label{lftadj} Suppose that the finite group $G$ acts on the $k$-linear category $\cc$. 
With the above notations one has that $\ind_H^G:
\cc^{H }\ra \cc^{G}$ is a $k$-linear left adjoint of the restriction functor $\Res^{G}_{H}:\cc^{G}\ra \cc^{H}$.\ep
\bpf
The proof of  \cite[Proposition 2.9]{buna} works verbatim in the case of an abelian category $\cc$.
\epf

\noindent
Note that the functor $\ind_H^G$ as defined above, depends on the set of representative elements  $\mtc R$. Since the adjoint of a functor is unique up to isomorphism it follows that for a different set of representative elements  one obtains an isomorphic functor. 
\subsection{Action by monoidal autoequivalences.}
Suppose that $\C$ is a $k$-linear monoidal category  and consider $\underline {\mtr{Aut}}_{\otimes}(\cc) $ the full subcategory of $\underline {Aut}(\cc)$ consisting of $k$-linear monoidal autoequivalences of $\cc$.

\noindent
We say that $G$ acts on $\C$ by \emph{monoidal autoequivalences}, if there is a unitary monoidal functor
$T: \underline G \to \underline {\mtr{Aut}}_{\otimes} \C$, Thus for any $g \in G$ there is,$T^{g}:\cc \ra \cc$ a monoidal autoequivalence of $\cc$ satisfying the conditions from above. Moreover, $T^g$ is endowed with a monoidal structure
$(T_2^g)^{M, N}: T^g(M \otimes N) \to T^g(M) \otimes
T^g(N)$, for all $M, N \in \C$ and $T_2^{g, h}: T^g T^h \to T^{gh}$ are natural isomorphisms of monoidal functors, for all $g, h \in G$. Thus, for
all $g, h\in G$ and $M, N \in \C$, Equation \eqref{tenant} becomes:
\begin{equation}\label{tensor-rho} {(T_2^{gh}})^{M, N} {(T_2^{g, h})}_{M \otimes N} =  ({(T_2^{g, h})}_{M} \otimes {(T_2^{g, h})}_{N}) \, {(T_2^{g})}^{T^h(M), T^h(N)} \, {T^g((T_2^{h})}_{M, N}).
\end{equation} 
\section{Categorical Mackey and Green functors}\label{catdef}
In this section we introduce the notion of categorical Mackey and Green functors and we prove some of their properties. First we recall the notion of Mackey and Green functors over rings.
\subsection{Classical Mackey functors} Let $G$ be a finite group.
A {\it Mackey functor} (or a $G$-functor) over a ring $R$ is a collection of $R$-modules $\{a(H)\}_{H\leq G}$
together 
with morphisms
$I^{H}_{K} : a(K)\ra a(H)$,
$R^{H}_{K}: a(H)\ra a(K)$ and
$c_{H, g} : a(H)\ra a( \;^{g}H)$
for all subgroups $H$ and $K$ of $G$ with $K\leq H$ and for all $g \in G$. This datum satisfies the following compatibility conditions:

$\mathrm{(M0)}$ $I_{H}^{H}, \;R_{H}^{H}, c_{H, h} : M(H)\ra M(H)$
 are the identity morphisms for all subgroups H and $h \in H$.

$\mathrm{(M1)}$ $R_{J}^{K}R_{K}^{H}=R^{H}_{J}$ for all subgroups $J\leq K\leq H$.

$\mathrm{(M2)}$ $
I_{K}^{H}I_{J}^{K}=I^{H}_{J}$, for all subgroups $J\leq K\leq H$.

$\mathrm{(M3)}$
$c_{H, g}c_{\;^{g}H, h} = c_{H, gh} $, for all $H \leq G$ and $g, h \in G$.

$\mathrm{(M4)}$
For any subgroups $K, H \leq G$ the following Mackey relation is satisfied:
\beq\label{mackeyss}
R^{G}_{H}I^{G}_{K}=\bigoplus_{x \in H\backslash G\slash K}I^{H}_{\;^{x}K\cap H}R^{^{x}K}_{^{x}K\cap H}c_{K, x}.
\eeq
Moreover, a {\it Green functor} over a commutative ring $R$, is a $G$-functor $a$
such that for any subgroup $H$ of $G$ one has that 
$a(H)$ is an associative $R$-algebra with identity and satisfying the following:

$\mathrm{(G1)}$ $R_{K}^{H}\;\;
\text{and} \;\; c_{{H, g}}\;\;\text{are always unitary R-algebra homomorphisms,}
$

$\mathrm{(G2)}$ $
I_{K}^{H}(aR^{H}_{K}(b)) =I_{K}^{H}(a)b
$,

$\mathrm{(G3)}$ $
I_{K}^{H}(R^{H}_{K}(b)a) =b I_{K}^{H}(a)
$
 for all subgroups $K\leq H$ and all $a \in a(K)$ and $b \in a(H)$.
\br
Note that the compatibility conditions $\mtr {(G2)}$ (and respectively $\mtr {(G3)}$) can be expressed as the fact that the induction maps $I^{H}_{K}$ are morphisms of right (and respectively left) $a(H)$-modules.
\er
\subsection{Categorical Mackey functors}
Let $\cS$ be a given $k$-linear monoidal category. A {\it categorical   Mackey functor} (or a categorical $G$-functor) over $\cs$ is a collection of 
\bne
\item $\cs$-module categories $\{\ct(L)\}_{L\leq G}$
\item $\cs$-module functors
$\mathcal I^{H}_{K} : \ct(K)\ra \ct(H)$,
\item $\cs$-module functors $\mathcal{R}^{H}_{K}: \ct(H)\ra \ct(K)$ and
\item $\cs$-module functors $c_{L,\; g} : \ct(L)\ra \ct( \;^{g}L)$ \ene
for all subgroups $K\leq H$ of $G$ and for all $g \in G$.  Moreover, the following compatibilities conditions are satisfied:

\noindent
$\mathrm{(CM0)}$ $\mathcal  I_{H}^{H}, \;\mathcal  R_{H}^{H}, c_{H, h} : \ct(H)\ra \ct(H)$
 are $k$-linear isomorphic to the identity morphisms, for any $h \in H$.

\noindent
$\mathrm{(CM1)}$ There are natural transformations which are isomorphisms of module functors $${\bf R}_{J, H}^{K}:\mathcal  R^{H}_{J} \xra{\simeq } \mathcal  R_{J}^{K}\mathcal  R_{K}^{H} $$  for all $J\leq K\leq H$

\noindent
$\mathrm{(CM2)}$ There are natural transformations which are isomorphisms of module functors $${\bf I}_{J, H}^{K}: \mathcal  I^{H}_{J}\xra {\simeq }
\mathcal I_{K}^{H}\mathcal  I_{J}^{K}$$ for all $J\leq K\leq H$.

\noindent
$\mathrm{(CM3)}$ There are natural transformations which are isomorphisms of module functors $${\bf C}^{H}_{a,b}:  c_{H, ab}\xra{\simeq } c_{\;^{b}H, a}c_{H, b},$$ for all $H \leq G$ and $a, b \in G$.

\noindent
$\mathrm{(CM4)}$ For any subgroups $K\leq H$
there are isomorphisms of module functors
\beqn
{\bf CI}_{H, a}^{K}:c_{H, a}\mathcal{I}_{K}^{H}\xra{\simeq } \mathcal{I}_{\:^{a}K}^{\;^{a}H}c_{K,a}
\eeqn

\noindent
$\mathrm{(CM5)}$ For any subgroups $K\leq H$ there are  natural transformation which are isomorphisms of module functors
\beqn
{\bf CR}_{H, a}^{K}:c_{K,a}\mathcal{R}_{K}^{H}\xra{\simeq } \mathcal{R}_{\:^{a}K}^{\;^{a}H}c_{H,a}
\eeqn

\noindent
$\mathrm{(CM6)}$ For any two subgroups $L,K\leq H$ one has an isomorphism
\beqn\label{mackeyss}
{\bf RI}^{H}_{L,K}:\mathcal  R^{H}_{L}\circ \mathcal  I^{H}_{K}\xra{\simeq} \bigoplus_{x \in L\backslash H\slash K}\mathcal  I^{L}_{\;^{x}K\cap L}\circ \mathcal  R^{^{x}K}_{^{x}K\cap L}\circ c_{K, x}.
\eeqn
as $\cs$-module functors.

\noindent
Moreover for any tower $J\leq K\leq L \leq H$ of subgroups of $G$ and any $a,b, c\in G$ one has that
 we have the following coherence relations between these natural transformations:
 \beq\label{idS} 
{\bf I}_{L, H}^{H}=\id_{\mathcal{I}_{L}^{H}}\;\;\text{and}\;\;{\bf I}_{L, H}^{L}=\id_{\mathcal{I}_{L}^{H}}
\eeq
\beq\label{Sid}
{\bf R}_{L, H}^{H}=\id_{\mathcal{R}_{L}^{H}}\;\;\text{and}\;\;{\bf R}_{L, H}^{L}=\id_{\mathcal{R}_{L}^{H}}
\eeq
\beq\label{unul}
{\bf C}_{a,1}^{H}=\id_{c_{H,a}}={\bf C}^{H}_{1,a}
\eeq
{\tiny
\begin{center}
\begin{equation}\tag{Diagram R}
\begin{tikzpicture}
  \matrix (m) [matrix of math nodes, row sep=8em,column sep=8em,minimum width=5em]
  {  \mathcal{R}_{J}^{K} \mathcal{R}_{K}^{L} \mathcal{R}_{L}^{H} &  \mathcal{R}_{J}^{K}\mathcal{R}_{K}^{H} \\
 \mathcal{R}_{J}^{L} \mathcal{R}_{L}^{H} & \mathcal{R}_{J}^{H}.  \\
    };
  \path[-stealth]
  (m-1-2) edge node [above] {$\;_{\mathcal{R}_{J}^{K}}({\bf R}_{K, H}^{L})$} (m-1-1)
    (m-2-2) edge node [above] {${\bf R}_{J, H}^{L}$} (m-2-1)
         
               (m-2-1) edge node [left] {$({\bf R}_{J, L}^{K})_{\mathcal R_{L}^{H}}$} (m-1-1)
    (m-2-2) edge node [right] {${\bf R}_{J, H}^{K}$} (m-1-2)
   ;
\end{tikzpicture}
\end{equation}
\end{center}
}
{\tiny
\begin{center}
\begin{equation}\tag{Diagram I}
\begin{tikzpicture}
  \matrix (m) [matrix of math nodes, row sep=8em,column sep=8em,minimum width=5em]
  {  \mathcal{I}_{L}^{H}\mathcal{I}_{K}^{L}\mathcal{I}_{J}^{K} & \mathcal{I}_{K}^{H} \mathcal{I}_{J}^{K}\\
  \mathcal{I}_{L}^{H}\mathcal{I}_{J}^{L} & \mathcal{I}_{J}^{H}.  \\
    };
  \path[-stealth]
  (m-1-2) edge node [above] {$({\bf I}_{K, H}^{L})_{\mathcal{I}_{J}^{K}}$} (m-1-1)
    (m-2-2) edge node [above] {${\bf I}_{J, H}^{L}$} (m-2-1)
         
               (m-2-1) edge node [left] {$_{\mathcal I_{L}^{H}}({\bf I}_{J, L}^{K})$} (m-1-1)
    (m-2-2) edge node [right] {${\bf I}_{J, H}^{K}$} (m-1-2)
   ;
\end{tikzpicture}
\end{equation}
\end{center}
}
{\tiny
\begin{center}
\begin{equation}\tag{Diagram C}
\begin{tikzpicture}
  \matrix (m) [matrix of math nodes, row sep=8em,column sep=8em,minimum width=5em]
  { c_{H, abc}& c_{\;^{c}H, ab}c_{H, c} \\
    c_{\;^{bc}H,a}c_{H, bc} & c_{\;^{bc}H,a}c_{\;^{c}H,b}c_{H,c}\\
    };
  \path[-stealth]
  (m-1-1) edge node [above] {${\bf C}^{H}_{ab,c}$} (m-1-2)
    (m-2-1) edge node [above] {$\;_{c_{\;^{bc}H, a}}({\bf C}^{H}_{b,c})$} (m-2-2)
               (m-1-1) edge node [left] {${\bf C}^{H}_{a,bc}$} (m-2-1)
    (m-1-2) edge node [right] {$({\bf C}^{\;^{c}H}_{a,b})_{c_{H,c}}$} (m-2-2)
   ;
\end{tikzpicture}
\end{equation}
\end{center}
\noindent
}
{\tiny
\begin{center}
\begin{equation}\tag{Diagram RRC}
\begin{tikzpicture}
  \matrix (m) [matrix of math nodes, row sep=6em,column sep=5em,minimum width=3em]
  { c_{K, a}\mathcal R^{L}_{K }\mathcal R_{L}^{H}& c_{L, a}\mathcal R^{H}_{L } & \\ && \mathcal R^{\;^{a}H}_{\;^{a}L}c_{H,a}\\
  \mathcal R^{\;^{a}L}_{\;^{a}K}  c_{L, a}\mathcal R^{H}_{L} & \mathcal R^{\;^{a}L}_{\;^{a}K} \mathcal{R}^{\;^{a}H}_{\;^{a}L} c_{H,a}\\
    };
  \path[-stealth]
  (m-1-2) edge node [above] {$\;_{c_{H,a}}{\bf R}_{K,H}^{L}$} (m-1-1)
    (m-3-1) edge node [above] {$\;_{\mathcal R_{\;^{a}L}^{\;^{a}R}}{\bf CR}^{L}_{H, a}$} (m-3-2)
    (m-2-3) edge node [right] {$({\bf R}^{\;^{a}L}_{\;^{a}K,\;^{a}H})_{c_{K, a}}$} (m-3-2)
         
               (m-1-1) edge node [left] {$(\mathrm{\bf CR}^{K}_{L,a})_{\mathcal R^{H}_{L}}$} (m-3-1)
    (m-1-2) edge node [right] {${\bf CR}^{L}_{H, a}$} (m-2-3)
   ;
\end{tikzpicture}
\end{equation}
\end{center}
}
{\tiny
\begin{center}
\begin{equation}\tag{Diagram RCC}
\begin{tikzpicture}
  \matrix (m) [matrix of math nodes, row sep=6em,column sep=4em,minimum width=3em]
  { c_{L, ab}\mathcal R_{L}^{H}&  \mathcal R^{\;^{ab}H}_{\;^{ab}L}c_{H,ab} & \\ && \mathcal R^{\;^{ab}H}_{\;^{ab}L}c_{\;^{b}L,a}c_{L,b}\\
    c_{\;^{b}L,a}c_{L, b}\mathcal R^{H}_{L}& c_{\;^{b}L,a}\mathcal R^{\;^{b}H}_{\;^{b}L}c_{L,b}&\\
    };
  \path[-stealth]
  (m-1-1) edge node [above] {${\bf CR}_{L, H}^{ab}$} (m-1-2)
    (m-3-1) edge node [above] {$\;_{c_{\;^{b}L,a}}({\bf CR}^{L}_{H, b})$} (m-3-2)
    (m-3-2) edge node [right] {$({\bf CR}^{\;^{a}}_{\;^{b}L,\;^{b}H,})_{c_{L, b}}$} (m-2-3)
               (m-1-1) edge node [left] {$({\bf C}^{L}_{a,b})_{\mathcal R^{H}_{L}}$} (m-3-1)
    (m-1-2) edge node [right] {$_{\mathcal R^{\;^{ab}H}_{\;^{ab}L}}({\bf C}^{L}_{a,b})$} (m-2-3)
   ;
\end{tikzpicture}
\end{equation}
\end{center}
}
\vskip -1cm
{\tiny
\begin{center}
\begin{equation}\tag{Diagram IIC}
\begin{tikzpicture}
  \matrix (m) [matrix of math nodes, row sep=6em,column sep=4em,minimum width=3em]
  { c_{H, a}\mathcal I_{L}^{H}\mathcal I^{L}_{K }& c_{H, a}\mathcal I^{H}_{K } & \\ && \mathcal I^{\;^{a}H}_{\;^{a}K}c_{K,a}\\
  \mathcal I^{\;^{a}H}_{\;^{a}L}  c_{L, a}\mathcal I^{L}_{K} & I^{\;^{a}H}_{\;^{a}L} I^{\;^{a}L}_{\;^{a}K} c_{K,a}\\
    };
  \path[-stealth]
  (m-1-2) edge node [above] {$\;_{c_{H,a}}({\bf R}_{K,H}^{L})$} (m-1-1)
    (m-3-1) edge node [above] {$\;_{\mathcal I_{\;^{a}L}^{\;^{a}H}}({\bf CI}^{K}_{L,a})$} (m-3-2)
    (m-2-3) edge node [right] {$({\bf R}^{\;^{a}L}_{\;^{a}K,\;^{a}H})_{c_{K, a}}$} (m-3-2)
         
               (m-1-1) edge node [left] {$({\bf CI}^{L}_{H,a})_{\mathcal I^{L}_{K}}$} (m-3-1)
    (m-1-2) edge node [right] {${\bf CI}^{K}_{H,a}$} (m-2-3)
   ;
\end{tikzpicture}
\end{equation}
   \end{center}
   }
{\tiny
\begin{center}
\begin{equation}\tag{Diagram ICC}
\begin{tikzpicture}
  \matrix (m) [matrix of math nodes, row sep=6em,column sep=4em,minimum width=3em]
  { c_{H, ab}\mathcal I_{L}^{H}&  \mathcal I^{\;^{ab}H}_{\;^{ab}L}c_{L,ab} & \\ && \mathcal I^{\;^{ab}H}_{\;^{ab}L}c_{\;^{b}L,a}c_{L,b}\\
    c_{\;^{b}H,a}c_{H, b}\mathcal I^{H}_{L}& c_{\;^{b}H,a}\mathcal I^{\;^{b}H}_{\;^{b}L}c_{L,b}&\\
    };
  \path[-stealth]
  (m-1-1) edge node [above] {${\bf CI}_{H, ab}^{L}$} (m-1-2)
    (m-3-1) edge node [above] {$\;_{c_{\;^{b}H,a}}({\bf CI}^{L}_{H,b})$} (m-3-2)
    (m-3-2) edge node [right] {$({\bf CI}^{\;^{b}L}_{\;^{b}H,a})_{c_{L, b}}$} (m-2-3)
         
               (m-1-1) edge node [left] {$({\bf C}^{H}_{a,b})_{\mathcal I^{H}_{L}}$} (m-3-1)
    (m-1-2) edge node [right] {$\;_{\mathcal I^{\;^{ab}H}_{\;^{ab}L}}({\bf C}^{L}_{a,b})$} (m-2-3)
   ;
\end{tikzpicture}
\end{equation}
\end{center}
}
\subsection{Categorical Green functors}
A {\it categorical Green functor} over $\cs$, is a Mackey-functor $\ct$ such that

\noindent $\mathrm{(CG0)}$
$\ct(L)$ is a unitary monoidal category,

\noindent
$\mathrm{(CG1)}$ $\mathcal R_{L}^{K}\;\;
\text{and} \;\; c_{{L, \;g}}\;\;\text{are always unitary monoidal functors,}
$

\noindent
$\mathrm{(CG2)}$ $
\mathcal I_{L}^{K}
$ is a $\ct(K)$-left module functor,

\noindent
$\mathrm{(CG3)}$ $
\mathcal I_{L}^{K}
$ is a $\ct(K)$-right module functor,

\noindent
$\mathrm{(CG4)}$ For any two subgroups $L,K\leq H$ one has that
\beq\label{mackeyssgr}
{\bf RI}^{H}_{L,K}:\mathcal  R^{H}_{L}\circ \mathcal  I^{H}_{K}\xra{\simeq} \bigoplus_{x \in L\backslash H\slash K}\mathcal  I^{L}_{\;^{x}K\cap L}\circ \mathcal  R^{^{x}K}_{^{x}K\cap L}\circ c_{K, x}.
\eeq
as $\ct(G)$-module functors.

\noindent
$\mathrm{(CG5)}$ The natural transformations ${\bf R}^{L}_{K, H}$, ${\bf C}^{H}_{a,b}$ are natural transformation of monoidal functors.

\noindent
$\mathrm{(CG6)}$ ${\bf I^{K}_{J,H}}$  is a  natural morphism $\ct(H)$ -module functors.

\noindent
$\mathrm{(CG7)}$ and ${\bf CI}^{K}_{H,a}, {\bf CR}^{K}_{H, a}$ are  natural morphisms of $\ct(G)$- module functors.
\section{Categorical Mackey functors from group actions on abelian categories}\label{tim}
The goal of this section is to prove the main results mentioned in the introduction.
\bp\lb{conj} 
Let $G$ be a finite group acting on the $k$-linear category $\cc$. Let $H\leq G$ be any subgroup of $G$ and $x \in G$.

\noindent
1) There is a $k$-linear functor $c_{H,x}:\cc^{H}\ra \cc^{\;^{x}H}$ which is a $k$-linear equivalence of categories.

\noindent 2) If $M=(V, \mu) \in \cc^H$ then $c_{H,x}(M):=(T^{x}(V), \;^{x}\mu)\in \cc^{\;^xH}$ with the equivariant structure $[\:^{x}\mu]^{xhx\inv}_{T^{x}(V)}:T^{xhx^{-1}}(T^x(V))\ra  T^x(V)$ given by
\beq\label{conjmu}
T^{xhx^{-1}}(T^x(V))\xra{(T_2^{xhx^{-1}, \;x})_{V}}T^{xh}(V)\xra{(T_2^{x,h})^{-1}_{V}}T^x(T^h(V))\xra{T^x(\mu^{h}_{V})} T^x(V),
\eeq
for any $h \in H$.

\noindent
 3) If $\cc$ is a $k$-linear monoidal category  and the action of $G$ on $\cc$ is by monoidal autoequivalences then $c_{H,x}$ is a $k$-linear monoidal functor. 
\ep
\bpf 1-2)
In order to see that $c_{H,x}(M)$ is an $\;^{x}H$-equivariant object it is enough to verify Equation \eqref{deltau} for any $xhx^{-1}, xlx^{-1}\in \;^{x}H$. This is equivalent to the diagram  below (made of solid arrows) being commutative. Note that compatibility conditions (\ref{ro-2})-(\ref{deltau}) of the action of $G$ imply the commutativity of the diagram after inserting the dashed arrows. Indeed, the bottom right trapeze $(5)$ is commutative by applying $T^{x}$ to the equivariantized condition (\ref{ro-2}) for $V\in \cc^{H}$. The adjacent trapeze $(6)$ is commutative by the naturallity of  $T_{2}^{{x,h}} $ with respect to the morphism $\mu_{V}^{l}$. The rectangle $(4)$ is commutative due to the associativity of the action, Equation \eqref{ro-2}. The parallelogram $(2)$ is commutative due to the associativity of the action, Equation \eqref{ro-2}. Diagram $(3)$ is commutative due to the  naturallity of the natural transformation $T^{xhx^{-1}, x}_{2}$ with respect to the morphism $T^{l}(V)\xra{\mu_{V}^{l}} V$. Diagram $(1)$ is commutative due to the associativity of the action,  Equation \eqref{ro-2}.
{\Tiny
 \begin{center}
{\Tiny
\begin{tikzpicture}
  \matrix (m) [matrix of math nodes, row sep=6em,column sep=1.15em,minimum width=3em]
  {
     T^{xhx\inv}(T^{xlx\inv}(T^{x}(V))) &  & T^{xhx\inv}(T^{xl}(V)) &  & T^{xhx\inv}(T^{x}(T^{l}(V))) \\
     & (1) & & (2) & \\
    T^{xhlx\inv} (T^{x}(V)) & &  & &  T^{xhx\inv}(T^{x}(V))\\
     & &  &  (3) & \\
     T^{xhl}(V)& & T^{xh}(T^{l}(V)) &  & T^{xh}(V)\\
     & (4) & & (6) &  \\
    T^{x}(T^{hl}(V)) & & T^{x}(T^{h}(T^{l}(V))) & & \\
    & (5) & & & \\
     T^{x}(V) & & & &  T^{x}(T^{h}(V)) \\
     };
  \path[-stealth]
    (m-1-1) edge node [left] {$(T_{2}^{xhx\inv, xlx\inv})_{V}$} (m-3-1)
    (m-1-3) edge [dashed] node [right] {$(T_{2}^{xhx\inv, xl})_{V}$} (m-5-1)
    (m-1-1) edge node [above] {$T^{xhx\inv}(T_{2}^{xlx\inv , x})_{V}$}(m-1-3)
  (m-3-5) edge node [right] {$(T_{2}^{xhx\inv, x})_{V}$} (m-5-5)
  (m-1-3) edge node [right] [above]{$T^{xhx\inv}(T_{2}^{x, l})_{V}$} (m-1-5)
         (m-1-5) edge node [right] {$T^{xhx\inv}(T^{x}(\mu_{V}^{l}))$} (m-3-5)
          (m-1-5) edge [dashed] node [left] {$(T_{2}^{xhx\inv, x})_{T^{l}(V)}$}
          (m-5-3)
           (m-5-3) edge [dashed] node [below] {$T^{xh}(\mu^{V}_{l})$}(m-5-5)
            (m-5-3) edge [dashed] node [right] {$(T_{2}^{x,h})_{T^{l}(V)}$}(m-7-3)
            (m-7-1) edge [dashed] node [below] {$T^{x}((T_{2}^{h,l})_{V})$} (m-7-3)
          (m-5-1) edge [dashed] node [below] {($T_{2}^{xh, l})^{-1}_{V}$} (m-5-3)
    (m-3-1) edge node [left] {$(T_{2}^{xhlx\inv , x})_{V}$} (m-5-1)
    (m-5-5) edge node [right] {$(T_{2}^{x,h})_{V}$} (m-9-5)
   (m-5-1) edge node [right] {$(T_{2}^{x, hl})^{-1}_{V}$} (m-7-1)
  (m-7-3) edge [dashed] node [right] {$T^{x}(T^{h}(\mu^{V}_{l}))$} (m-9-5)
   (m-7-1) edge node [right] {$T^{x}(\mu^{hl}_{V})$} (m-9-1)
      (m-9-5) edge node [below] {$T^{x}(\mu^{h}_{V})$} (m-9-1)
         ;
\end{tikzpicture}
}
\end{center}
}
\noindent It is easy also to verify that if $f:M \ra N$ is  a morphism in $\cc^{H}$ then $T^{x}(f)$ is a morphism in $\cc^{\;^{x}H}$.
Then it is clear that $c_{H,x}$ is an equivalence of categories with the inverse given by $c_{\;^{x}H,\;x^{-1}}:\cc^{\;^{x}H}\ra \cc^{H}$.

\noindent 2) Let $M, N$ be two objects of $\cc^{H}$. If the action of $G$ on $\cc$ is by monoidal autoequivalences then one can consider $(T^{x}_{2})^{M,N}$ as the monoidal structure of $c_{H,x}:c_{H, x}(M\ot N)\ra c_{H,x}(M)\ot c_{H,x}(N)$. It is clear that  $(T^{x}_{2})^{M,N}$ satisfies the pentagon axiom from the definition of a monoidal structure.

\noindent
One has to check that the monoidal structure $(T^{x}_{2})^{M, N} :T^{x}(M\ot N)\ra T^{x}(M)\ot T^{x}(N)$ of $T^{x}$ is a morphism in $\cc^{\;^{xH}}$. Thus for any $M, N \in \cc^{\;^{x}H}$ one has to check the commutativity of the following digram:

\begin{center}
{\Tiny
\begin{tikzpicture}
  \matrix (m) [matrix of math nodes, row sep=2em,column sep= 0.9 em,minimum width=3em]
  { T^{xhx\inv}(T^{x}(M\ot N)) &   &  T^{xhx\inv}(T^{x}(M)\ot T^{x}(N))\\
  & (1) & \\ 
   &  & T^{xhx\inv}(T^{x}(M))\ot T^{xhx\inv}(T^{x}(N))\\
   T^{xh}(M\ot N) & & T^{xh}(M)\ot T^{xh}(N)\\
    & (2) & \\
  T^{x}(T^{h}(M\ot N)) & & \\
  T^{x}(T^{h}(M)\ot T^{h}(N)) & & T^{x}(T^{h}(M))\ot T^{x}(T^{h}(N))  \\
  & (3) & \\
  T^{x}(M\ot N) & &T^{x}(M)\ot T^{x}(N)\\};
  \path[-stealth]
  (m-1-1) edge node [above] {$T^{xhx\inv}((T^{x}_{2})^{M, N})$} (m-1-3)
    (m-9-1) edge node [above] {${(T^{x}_{2})}^{M, N}$} (m-9-3)
         (m-4-1) edge node [left] {$(T^{x,h}_{2})^{-1}_{M\ot N}$} (m-6-1)
         (m-6-1) edge node [left] {$T^{x}((T_{2}^{h})^{M, N})$} (m-7-1)
         (m-7-1) edge node [left] {$T^{x}(\mu^{M}_{h}\ot \mu^{N}_{h})$} (m-9-1)
            (m-3-3) edge node [right] {$(T_{2}^{xhx\inv, x})_{M}\ot (T_{2}^{xhx\inv, x})_{N}$} (m-4-3)
         (m-4-3) edge node [right] {$(T_{2}^{x, h})^{{-1}}_{M}\ot (T_{2}^{x, h})^{{-1}}_{N}$} (m-7-3)
         (m-7-3) edge node [right] {$ T^{x}(\mu^{M}_{h}) \ot T^{x}(\mu^{N}_{h})$} (m-9-3)
               (m-1-1) edge node [left] {$(T^{xhx\inv, x}_{2})_{M\ot N}$} (m-4-1)
    (m-1-3) edge node [right] {$(T^{xhx\inv}_{2})^{{T^{x}(M), T^{x}(N)}}$} (m-3-3)
    (m-4-1) edge [dashed] node [above] {$(T^{xh}_{2})^{M,N}$} (m-4-3)
     (m-7-1) edge [dashed] node [above] {$({T_{2}^{x})}^{{T^{h}(M)},{T^{h}(N)}}$} (m-7-3)
    ;
\end{tikzpicture}
}
\end{center}
The upper pentagon (1) is commutative since $T_{2}^{xax\inv , \;x}$ is a natural transformation of monoidal functors, Equation \eqref{tensor-rho}. The middle pentagon $(2)$ is commutative since $T^{x, a}_{2}$ is a natural transformation of monoidal functors, same Equation \eqref{tensor-rho}. The bottom rectangle commutes from the compatibility condition of the monoidal functor $T^{x}$ with the tensor product of morphisms, Equation \eqref{natf}. \epf
\br\label{idc}
Note that if $x \in H$ then $c_{H, x }\simeq \id_{\cc^{H}}$ via the natural transformation $N_{V}:c_{H, x}(V) \xra{\mu_{V}^{x}} V$ for any $(V, \mu)\in \cc^{H}$. Moreover, if $\cc$ is a $k$-linear monoidal category  and the action of $G$ on $\cc$ is by monoidal autoequivalences then $c_{H,x} \simeq \id_{\cc^{H}}$ as monoidal functors. 
\er
\bp\label{nat-c} Let $G$ be a finite group acting $k$-linearly on a $k$-linear category $\cc$.
With the above notations there is a natural transformation which is an  isomorphism of $k$-linear functors\beqn
{\bf C}^{H}_{a,b}:c_{H, ab}\ra c_{\;^{b}H, a}c_{H, b}
\eeqn
\noindent
 If $\cc$ is a $k$-linear monoidal category  and the action of $G$ on $\cc$ is by monoidal autoequivalences then ${\bf C}^{H}_{a,b}$ is an isomorphism of monoidal functors.
\ep
\bpf
One can easily check that $({\bf C}^{H}_{a,b})_{M}:=(T^{a,b}_{2})^{-1}_{M}\colon T^{ab}(M)\ra T^{a}(T^{b}(M))$ is a morphism in $\cc^{\;^{ab}H}$. Indeed for any $h \in H$ one has to check that the following diagram is commutative:
\begin{center}
{\Tiny
\begin{tikzpicture}
  \matrix (m) [matrix of math nodes, row sep=5em,column sep=7em,minimum width=4em]
  {T^{abh(ab)\inv}(T^{ab}(M) ) &  T^{abh(ab)\inv}(T^{a}(T^{b}(M) )) \\
T^{abh}(M))  &   T^{abhb^{-1}}(T^{b}(M))\\ 
   T^{ab}(T^{h}(M)) &  T^{a}(T^{bhb^{-1}}(T^{b}(M)))\\
  T^{ab}(M) &  T^{a}(T^{bh}(M))\\
  &  T^{a}(T^{b}(T^{h}(M))) \\
  &  T^{a}(T^{b}(M)) \\
};
  \path[-stealth] 
     (m-1-1) edge  node [above] {$T^{abh(ab)\inv}((T_{2}^{a, b})^{-1}_{M})$} (m-1-2)
  (m-1-1) edge node [left] {$(T_{2}^{abh(ab)^{-1},ab})_{M}$} (m-2-1) 
  (m-2-1) edge node [left] {$(T_{2}^{ab,h})^{-1}_{M}$} (m-3-1)
   (m-3-1) edge node [left] {$T^{ab}(\mu_{M}^{h})$} (m-4-1)
      (m-4-1) edge node [below] {$(T_{2}^{a, b})^{-1}_{M}$} (m-6-2)
  (m-1-2) edge  node [right] {$(T^{abh(ab)^{-1}, a}_{2})_{T^{b}(M)}$} (m-2-2)
  (m-2-2) edge  node [right] {$(T^{abhb^{-1}, b}_{2})_{M}$} (m-3-2)
  (m-3-2) edge node [right] {$T^{a}((T^{bhb\inv, b}_{2})_{M})$} (m-4-2)
   (m-4-2) edge  node [right] {$T^{a}((T^{b, h}_{2})^{-1}_{M})$} (m-5-2)
    (m-2-1) edge [dashed] node [above] {$$} (m-2-2)
     (m-2-1) edge [dashed] node [above] {$$} (m-4-2)
     (m-3-1) edge [dashed] node [above] {$$} (m-5-2)
  (m-5-2) edge node [right] {$T^{a}(T^{b}(\mu_{M}^{h}))$} (m-6-2)
    ;
\end{tikzpicture}
}
\end{center}
\noindent
Note that the last parallelogram is commutative by the naturality of $T^{a, b}_{2}$. All the other parallelograms are commutative by Equation \eqref{deltau}. This defines a natural transformation ${\bf C}^{H}_{a,b}:c_{H, ab}\ra c_{\;^{b}H, a}c_{H, b}$ which clearly it is an isomorphism.

\medbreak\noindent
If the action is by monoidal equivalences then it is clear that ${\bf C}^{H}_{a,b}$ is a  natural morphisms of monoidal functors. Indeed $C^{H}_{a,b}$ satisfies Equation \eqref{tenant} since $T^{a,b}_{2}:T^{a}T^{b}\ra T^{ab}$ is also a morphism of monoidal functors.
\epf
\subsection{Module category structures} Suppose that the group $G$ acts by monoidal equivalences on the $k$-linear monoidal category  $\cc$. For any $L \leq H \leq G$ note that $\cc^{L}$ is a $\cc^{H}$-left (or right) module category via the restriction functor $\Res_{H}^{L}$. Thus one can define 
$M\boxtimes (V ,\mu):=(\Res_{H}^{L}(M)\ot V, \nu)$ where for any $h \in H$ one has that $\nu^{h}_{M\ot V}:T^{h}(M\ot V) \ra M\ot V$ is defined via: $$\nu^{h}_{M\ot V}: T^{h}(M\otimes V) \xra{(T^{h}_{2})^{M,V}} T^{h}(M)\ot T^{h}(V)\xra{\mu_{V}^{h}\otimes \mu_{M}^{h}} M \otimes V.$$
Note also that $M\boxtimes V=\res^{G}_{H}(M)\ot V$ 
and  the $\cc^{L}$-module category associativity constraint of $\cc^{H}$ coincides to the associativity constraint of $\cc$ as a monoidal category. 
\subsubsection{On the module functor structures of the restriction and induction functors} 
If $L\leq H$ are subgroups of $G$ then the restriction functor $\res^{H}_{L}:\cc^{H}\ra \cc^{L}$ is a $\cc^{G}$-module  functor. Indeed note that $\re^{H}_{L}(V\boxtimes M)=\re^{H}_{L}(\re^{G}_{H}(M)\ot V)=\re^{H}_{L}(\re^{G}_{H})(M)\ot \re^{H}_{L}(V)=\re^{G}_{L}(M)\ot \re^{H}_{L}(V)=M\boxtimes \re^{H}_{L}(V)$.

\noindent  Thus the module functor structure of the functor $\re^{H}_{L}$ can be considered as the identity map. Moreover,  since for any $K\leq L\leq H$ one has $\re^{L}_{K}\re^{H}_{L}=\re^{H}_{K}$ it follows that one can also consider ${\bf R}^{L}_{K, H}:\re^{L}_{K}\re^{H}_{L}\ra \re^{H}_{K}$ as the identity natural transformation.

\noindent
Let $(V, \mu) \in \cc^{H}$ and $L\leq H$. By abuse of notations sometimes we write $\ind_{H}^{G}(V)$ instead of $\ind_{H}^{G}(V, \mu)$. We also write shortly $\re^{H}_{L}(V)$ instead of $\re^{H}_{L}(V, \mu)$.
\bl Let $G$ be a finite group acting by monoidal autoequivalences on $\cc$ and $L\leq H\leq G$ be a tower of subgroups. With the above notations  the induction functor $\ind_{L}^{H}:\cc^{L}\ra \cc^{H}$ is a $\cc^{H}$-module functor.
\el 
\bpf 
Let $M \in \cc^{H}$ and $V \in \cc^{L}$. We will define a canonical isomorphism
\beq\label{is}
(\ind^{H}_{L})^{M, V}_{2}:\ind^{H}_{L}( M\boxtimes V)\ra M\boxtimes \ind^{H}_{L}(V)
\eeq
in $\cc^{H}$ and we will show that it satisfies the module functor axioms.  Fix a set $\mtc R$ of representative elements for the left cosets of $H\leq L$. Let $\ind^{L}_{H}$ be the associated induction functor $\ind_{L}^{H}:\cc^{L}\ra \cc^{H}$ to the set $\mtc R$.  Thus
\beqarn
\ind^{L}_{H}(M\boxtimes V)=\ind^{L}_{H}(\re^{L}_{H}(M)\ot V)=\oplus_{a \in \mathcal R}T^{a}(M\otimes V)
\eeqarn
and 
\beqn
M\boxtimes \ind^{H}_{L}(V)= \re^{G}_{H}(M)\ot (\oplus_{a \in H/L}T^{a}(V))=\oplus_{a \in H/L} M\ot T^{a}(V)
\eeqn 
On the components, for any $a \in H$ the module functor structure will be defined as \beqn
{(\ind^{H}_{L})^{M, V}_{2}}^{\;,a}:
T^{a}(M\otimes V)\xra{(T^{a}_{2})^{M,V}} T^{a}(M)\ot T^{a}(V)
 \xra{\mu_{M}^{a}\ot 1}M\ot T^{a}(V).
 \eeqn
One needs to verify that the above module functor structure is a morphism in $\cc^{H}$ which is equivalent to the following commutative diagram:
{\tiny
\begin{center}
\begin{tikzpicture}
  \matrix (m) [matrix of math nodes, row sep=9em,column sep=6em,minimum width=5em]
  { 
T^{h}(\ind^{H}_{L}(M\boxtimes V) &
 T^{h}(M\boxtimes \ind^{H}_{L}(V))\\
    \ind^{H}_{L}(M\boxtimes V) & M\boxtimes \ind^{H}_{L}(V).\\
    };
  \path[-stealth]
  (m-1-1) edge node [above] {$T^{h}((\ind^{H}_{L})^{M, V}_{2})$} (m-1-2)
    (m-2-1) edge node [above] {$(\ind^{H}_{L})^{M, V}_{2}$} (m-2-2)
         
               (m-1-1) edge node [left] {$\mu^{h}_{\ind^{H}_{L}(M\boxtimes V)}$} (m-2-1)
    (m-1-2) edge node [right] {$\mu^{h}_{M\boxtimes \ind^{H}_{L}(V)}$} (m-2-2)
   ;
\end{tikzpicture}
   \end{center}
}
On the components this can be written as the commutativity of the following:
\begin{center}
{\Tiny
\begin{tikzpicture}
  \matrix (m) [matrix of math nodes, row sep=3.25 em,column sep=0.00005 em,minimum width=0.1em]
  { 
     T^{h}(T^{a}(M\ot V)) &  & T^{h}(T^{a}(M)\ot T^{a}(V)) &  & T^{h}(M\ot T^{a}(V)) \\
     & (1) & T^{h}(T^{a}(M))\ot T^{h}(T^{a}(V))&\; & \\
    T^{ha} (M\ot  V) & & T^{ha}(M) \ot T^{ha}(V)  &  & T^{h}(M)\ot T^{h}(T^{a}(V))\\
     & (3) & & (4) & \\
     T^{b}(T^{l}(M\ot V))& &  &  & M\ot T^{h}(T^{a}(V))\\
    T^{b}(T^{l}(M)\ot T^{l}(V)) &  & T^{b}(T^{l}(M))\ot T^{b}(T^{l}(V)) & (6)   & M\ot T^{ha}(V) \\
     T^{b} (M\ot  V) & & T^{b}(M)\ot T^{b}(T^{l}(V)) & & M\ot T^{b}(T^{l}(V)) \\
      T^{b}(M)\ot T^{b}(V) & &  & &  M\ot  T^{b}(V) \\
     };
  \path[-stealth]
    (m-1-1) edge node [left] {$(T_{2}^{h, x})_{T^{r}(M)}$} (m-3-1)
    (m-3-3) edge [dashed] node [left] {$(T_{2}^{b, l})^{-1}_{M}\ot (T_{2}^{b, l})^{-1}_{V}$} (m-6-3)
    (m-1-1) edge node [above] {$T^{h}[(T_{2}^{a})^{M, V}]$}(m-1-3)
  (m-3-5) edge node [right] {$\mu_{M}^{h}1$} (m-5-5)
  (m-1-3) edge node [right] [above]{${T^{h}(\mu^{a}_{M}1)}$} (m-1-5)
         (m-1-5) edge node [right] {$(T_{2}^{h})^{M, T^{a}(V)}$} (m-3-5)
          (m-6-1) edge [dashed] node [above] {$(T_{2}^{b})_{T^{l}(M),T^{l}(V)}$}(m-6-3)
           (m-7-3) edge [dashed] node [below] {$1\ot T^{b}(\mu^{l}_{V})$}(m-8-1)
            (m-7-3) edge [dashed] node [below] {$\mu_{M}^{b}\ot 1$}(m-7-5)
            (m-2-3) edge [dashed] node [below] {$T^{h}(\mu_{M}^{a})1$} (m-3-5)
             (m-1-3) edge [dashed] node [right] {$(T_{2}^{h})^{T^{a}(M), T^{a}(V)}$} (m-2-3)
             (m-2-3) edge [dashed] node [left] {$(T_{2}^{h, a})_{M}(T_{2}^{h, a})_{V}$} (m-3-3)
         (m-6-3) edge [dashed] node [left] {$$} (m-7-3)
    (m-3-1) edge node [left] {$(T_{2}^{x, a})^{-1}_{T^{r}(M)}$} (m-5-1)
      (m-7-5) edge node [right] {$1T^{b}(\mu^{l}_{V}))$} (m-8-5)
(m-6-5) edge node [left] {$M((T_{2}^{b, l})^{-1}_{V})$} (m-7-5)   
(m-5-5) edge node [right] {$1(T_{2}^{h, a})_{V}$} (m-6-5)
(m-5-1) edge node [left] {$T^{b}((T_{2}^{l})^{M, V})$} (m-6-1)
(m-6-1) edge node [left] {$T^{b}(\mu^{l}_{M}\ot \mu^{l}_{V})$} (m-7-1)
    (m-8-1) edge node [below] {$\mu_{M}^{b}\ot 1$} (m-8-5)
(m-3-1) edge [dashed] node [below] {$(T^{ha}_{2})^{M, V}$} (m-3-3)
(m-3-3) edge [dashed] node [left] [above] {$\mu_{M}^{ha}1$} (m-6-5)
   (m-7-1) edge node [left] {$(T^{b}_{2})^{M, V}$} (m-8-1)
         ;
\end{tikzpicture}
}
\end{center}
\noindent Using Equation \eqref{deltau} it is straightforward to verify that the above map $(\ind^{H}_{L})^{M, V}_{2}$ is a module functor structure.
\epf
\br\label{remf}
Note that the induction functor defined above depends on the set $\mtc R$ of chosen representative elements  for the left cosets $L \leq H$. Moreover, as explained in Section \ref{indc} changing the representative elements  one obtains a $k$-linear isomorphic functor. Since $\Res^{H}_{L}$ is a $\cc^{H}$-module functor and $\ind^{L}_{H}$ is its left adjoint it follows that $\ind^{H}_{L}$ {\it as a $\cc^{H}$-module functor} does not depend on  the chosen set $\mtc R$. \er

\bl\lb{commir} Let $G$ be a finite group acting on the abelian category $\cc$.

\noindent
1) Suppose that $L\leq H$ are two subgroups of $G$. Then the following identities hold for any $x \in G$:
\beq\label{cr}
c_{L, x}\circ \res_{L}^H=\res_{\;^xL}^{\;^xH} \circ c_{H, x}
\eq	
as functors from $\cc^{L}$ to $\cc^{\;^{x}H}$.

\noindent 
2) There is a natural transformation 
\bq
{\bf CI}_{L, H}^{x}: c_{H, x}\circ \ind_{L}^H\xra{\simeq} \ind_{\;^xL}^{\;^xH} \circ c_{L, x}
\eq 
which is a natural isomorphism of functors from $\cc^{L}$ to $\cc^{\;^{x}H}$

\noindent
3) If $\cc$ is a $k$-linear monoidal category  and the action of $G$ on $\cc$ is by monoidal autoequivalences then ${\bf CI}_{L, H}^{x}$ is an isomorphism of $\cc^{G}$-module functors.
\el
\bpf
1) The first identity is straightforward. 

\noindent
2) Let $\mtc{R}$ be a set of representative elements for the left cosets $\{Lx\;|x \in H\}$ of the extension $H/L$. Suppose that 
\bq
\ind_L^H(M)=(\opl_{r \in \mtc{R}}T^r(M), \nu_{\ind_L^H(M)})
\eq
for any  $(M, \{\mu^{b}_{M}\}_{b\in L})\in \cc^L$. Then using formula (\ref{muind}) it follows that $\nu_{\ind_L^H(M)}^{a}$ is defined on the components as follows:
\bq
\nu^{a, r}:T^a(T^r(M))\xra{T^{a,r}_{2}(M)} T^{ar}(M)\xra{{(T_{2}^{r', b}})_{M}^{-1}} T^{r'}(T^b(M))\xra{T^{r'}(\mu_M^{b})
} T^{r'}(M)\eq
where $r' \in \mtc{R}$ and $b \in L$ are determined by $ar=r'b$. On the other hand, using Proposition \ref{conj} for $c_{H, x}(\ind^{H}_{L}(M))$, as an object of $\cc^{\;^{x}H}$ one has that
\bq
c_{H, x}(\ind^{H}_{L}(M))=(\oplus_{r \in \mtc{R}}T^{x}(T^{r}(M)), \;^{x}\nu_{T^{x}(\ind^{H}_{L}(M))})
\eq
with the equivariant structure $[\;^{x}\nu]$ given on the components  by:
\beqn
 [\;^{x}\nu]^{xax^{-1}, r}:T^{xax^{-1}}(T^{x}(T^{r}(M))) \xra{({T_{2}^{{h, x}} })_{T^{r}(M)}} T^{xa}(T^{r}(M))\xra{{(T_{2}^{x,a})^{-1}}_{T^{r}(M)}}\eeqn \beqn  \xra{{(T_{2}^{x,a})^{-1}}_{T^{r}(M)}}  T^{x}(T^{a}(T^{r}(M))) \xra{T^{x}(T_{2}^{a, r})_{M}} T^{x}(T^{ar}(M))\xra{T^{x}(T_{2}^{r', b})^{{-1}}_{M}}\eeqn \beqn \xra{T^{x}(T_{2}^{r', b})^{{-1}}_{M}} T^{x}(T^{r'}(T^{b}(M)))\xra{T^{x}(T^{r'}(\mu^{b}_{M}))} T^{x}(T^{{r'}}(M))
\eeqn
\noindent
On the other hand, since $x\mtc{R}x\inv$ is a set of representative for the left cosets of $\;^{x}H/\;^{x}L$, by Remark \ref{remf} one may suppose  that 
\bq
\ind_{\;^xL}^{\;^xH}c_{L, x}(M)=(\oplus_{r \in \mtc{R}}T^{xrx\inv}(T^{x}(M)), \;\eta^{}_{\ind_{\;^xL}^{\;^xH}T^{x}(M)})
\eq
Using again formula (\ref{muind}) it follows that the equivariant structure of $\ind_{\;^xL}^{\;^xH}T^{x}(M)$ on the components is given as follows:
\beqn
\eta^{xax^{-1}, \;xrx^{-1}}_{M} :T^{xax^{-1}}(T^{xrx\inv}(T^{x}(M)) \xra{({T_{2}^{{xax\inv, xrx\inv}} })_{T^{x}(M)}} \eeqn \beqn T^{xarx\inv}(T^{x}(M)) \xra{{(T_{2}^{xr'x\inv,xbx\inv})^{-1}}_{T^{x}(M)}}  T^{xr'x\inv}(T^{xbx\inv }(T^{x}(M)))\eeqn \beqn \xra{T^{xr'x\inv}((\;^{x}\mu)^{xbx\inv})} T^{xr'x\inv}(T^{x}(M))
\eeqn
Define the natural transformation $ \mathbf{CI}^{x}_{L, H}:=\oplus_{r \in \mtc{R}}(T_{2}^{xrx\inv, x})^{-1}(T_{2}^{x,r}):c_{H, x}\ind_{L}^{H}\ra \ind_{\;^xL}^{\;^xH}c_{L,x}$. 
It will be shown that $({\mathbf{CI}^{x}_{L, H}})_{M}:T^{x}(\ind_{L}^{H}(M))\ra \ind_{\;^xL}^{\;^xH}(T^{x}(M))$ is an isomorphism in $\cc^{\;^{x}H}$ for any $M \in \cc^{L}$.  Indeed, one has to verify that the above morphism $\mathbf{CI}^{x}_{L, H}$ is compatible with the two equivariant structures defined above. This means that the following diagram:
is commutative.    {\tiny
\begin{center}
\begin{tikzpicture}
  \matrix (m) [matrix of math nodes, row sep=8em,column sep=6em,minimum width=5em]
  { T^{xax^{-1}}(T^{x}\ind^{H}_{L}(M)) & T^{xax^{-1}}( \ind^{\;^{x}H}_{\;^{x}L}T^{x}(M)) \\
    T^{x}(\ind^{H}_{L}(M)) & \ind^{\;^{x}H}_{\;^{x}L}T^{x}(M).\\
    };
  \path[-stealth]
  (m-1-1) edge node [above] {$T^{xax^{-1}}((\mathbf{CI}^{x}_{L, H})_{M})$} (m-1-2)
    (m-2-1) edge node [above] {$(\mathbf{CI}^{x}_{L, H})_{M}$} (m-2-2)
         (m-1-1) edge node [left] {$(\;^{x}\nu^{xax^{-1}})_{M}$} (m-2-1)
    (m-1-2) edge node [right] {$\eta^{xax^{-1}}_{M}$} (m-2-2)
   ;
\end{tikzpicture}
   \end{center}
}
On the components the above diagram becomes the following: \begin{center}
{\tiny
\begin{tikzpicture}
  \matrix (m) [matrix of math nodes, row sep=3.25 em,column sep=0.00005 em,minimum width=0.1em]
  { 
     T^{xax\inv}(T^{x}(T^{r}(M))) &  & T^{xax\inv}(T^{xr}(M)) &  & T^{xax\inv}(T^{xrx^{-1}}(T^{x}(M))) \\
     & (1) & &\;\; \;\;(2) & \\
    T^{xa} (T^{r}(M)) & & T^{xar}(M) =T^{xr'b}(M)  &  & T^{xarx\inv}(T^{x}(M))=T^{xr'bx\inv}(T^{x}(M))\\
     & (3) & & (4) & \\
     T^{x}(T^{a}(T^{r}(M)))& &  &  & T^{xr'x^{-1}}(T^{xbx^{-1}}(T^{x}(M)))\\
    T^{x}(T^{ar}(M))=T^{x}(T^{r'b}(M)) & (5) & & (6)   & T^{xr'x^{-1}}(T^{xb}(M)) \\
     T^{x}(T^{r'}(T^{b}(M)))) & & T^{xr'}(T^{b}(M)) & &  T^{xr'x^{-1}}(T^{x}(T^{b}(M))) \\
      T^{x}(T^{r'}(M))) & & T^{xr'}(M) & &  T^{xr'x^{-1}}(T^{x}(M)) \\
     };
  \path[-stealth]
    (m-1-1) edge node [left] {$(T_{2}^{xax\inv, x})_{T^{r}(M)}$} (m-3-1)
    (m-1-1) edge node [above] {$T^{xax\inv}(T_{2}^{x , r})_{M}$}(m-1-3)
  (m-3-5) edge node [right] {$(T_{2}^{xr'x\inv, xbx^{-1}})^{-1}_{T^{x}(M)}$} (m-5-5)
  (m-1-3) edge node [right] [above]{$(T_{2}^{x, a})^{-1}_{T^{r}(M)}$} (m-1-5)
         (m-1-5) edge node [right] {$(T_{2}^{xax\inv, xrx^{-1}})_{T^{x}(M)}$} (m-3-5)
          (m-5-3)
            (m-7-3) edge [dashed] node [right] {$(T_{2}^{xr' b})_{M}$}(m-3-3)
            (m-7-1) edge [dashed] node [below] {$(T_{2}^{x,r'})_{T^{b}(M)})$} (m-7-3)
             (m-1-3) edge [dashed] node [right] {$(T_{2}^{xax^{-1},xr})_{M}$} (m-3-3)
    (m-3-1) edge node [left] {$(T_{2}^{x, a})^{-1}_{T^{r}(M)}$} (m-5-1)
      (m-7-3) edge [dashed]  node [right] {$T^{xr'}(\mu^{b}_{M})$} (m-8-3)
      (m-7-5) edge node [right] {$T^{xr'x^{-1}}(T^{x}(\mu^{b}_{M}))$} (m-8-5)
(m-6-5) edge node [left] {$T^{xr'x^{-1}}((T_{2}^{x, b})^{-1}_{M})$} (m-7-5)   
(m-5-5) edge node [right] {$T^{xr'x^{-1}}((T_{2}^{xbx^{-1}, x})_{M})$} (m-6-5)
(m-5-1) edge node [left] {$T^{x}((T_{2}^{a, r})^{-1}_{M})$} (m-6-1)
(m-6-1) edge node [left] {$T^{x}((T_{2}^{r', b})^{-1}_{M})$} (m-7-1)
  (m-7-3) edge [dashed] node [below] {$(T^{xr'x^{-1}, x}_{2})^{-1}_{T^{b}(M)}$} (m-7-5)
    (m-8-1) edge node [below] {$(T^{x, r'}_{2})_{T^{b}(M)}$} (m-8-3)
(m-8-3) edge node [below] {$(T^{xr'x^{-1}, x}_{2})^{-1}_{M}$} (m-8-5)
(m-3-1) edge [dashed] node [below] {$(T^{xa, r}_{2})_{M}$} (m-3-3)
(m-3-3) edge [dashed] node [below] {$(T^{xr'x^{-1}, x}_{2})^{-1}_{T^{b}(M)}$} (m-3-5)
(m-3-3) edge [dashed] node [left] [above] {$(T^{xr'x^{-1}, xb}_{2})^{-1}_{M}$} (m-6-5)
(m-3-3) edge [dashed] node [below] {$(T^{x, r'b}_{2})^{-1}_{T^{b}(M)}$} (m-6-1)
   (m-7-1) edge node [left] {$T^{x}(T^{r'}(\mu^{b}_{M}))$} (m-8-1)
         ;
\end{tikzpicture}
}
\end{center}
Note that diagrams $(1)-(6)$ are commutativity by the associativity of the action, Equation \eqref{ro-2}. The bottom two rectangles are commutative since $T^{x, r'}_{2}$ and $T^{xr'x^{-1}, x}_{2}$ are natural transformations.
\vskip 0,5 cm
\noindent
 3) It is also straightforward to verify that this defines a natural transformation of $\cc^{G}$-module functors in the case of an action by monoidal autoequivalences. Indeed let $M \in \cc^{G}$ and $V \in \cc^{H}$. One has to check that the following diagram is commutative:
{
\begin{center}
{\tiny
\begin{tikzpicture}
  \matrix (m) [matrix of math nodes, row sep=7em,column sep=8em,minimum width=5em]
  { c_{L, x}\ind^{L}_{H}(M\ot V) & \ind^{\;^{x}L}_{\;^{x}H}c_{H, x}(M\ot V) \\
    M\ot c_{L, x}\ind^{L}_{H}( V) & M\ot \ind^{\;^{x}L}_{\;^{x}H}c_{H, x}( V)\\
    };
  \path[-stealth]
  (m-1-1) edge node [above] {$(\mathbf{CI}^{x}_{L, H})_{M\ot V}$} (m-1-2)
    (m-2-1) edge node [above] {$1\ot (\mathbf{CI}^{x}_{L, H})_{V}$} (m-2-2)
         (m-1-1) edge node [left] {$(c_{L, x}\ind^{L}_{H})_{2}^{M,V}$} (m-2-1)
    (m-1-2) edge node [right] {$(\ind^{\;^{x}L}_{\;^{x}H}c_{H, x})_{2}^{M,V}$} (m-2-2)
   ;
\end{tikzpicture}
}
   \end{center}
}

On the components the above diagram is equivalent to the commutativity of the following diagram; note that for shortness we replaced the symbol ``$\ot$'' by ``$.$''.
\begin{center}
{\Tiny
\begin{tikzpicture}
  \matrix (m) [matrix of math nodes, row sep=6.25 em,column sep=0.00015 em,minimum width=0.1em]
  { 
     T^{x}(T^{r}(M. V)) &  & T^{xr}(M. V) &   & & T^{xrx^{-1}}(T^{x}(M. V))\\
     & (D1) & T^{xr}(M) . T^{xr}(V) & (D2) &    &\\
    T^{x} (T^{r}(M). T^{r}(V)) & & T^{x}(T^{r}(M)) . T^{x}(T^{r}(V) )& &  & T^{xrx^{-1}}(T^{x}(M). T^{x}( V)) \\
     & (D3) &     & T^{xrx^{-1}}(T^{x}(M)). T^{xrx^{-1}}(T^{x}( V) )&   (D4)  & \\
     T^{x}(M. T^{r}(V)) & &   &  & & T^{xrx^{-1}}(M. T^{x}( V))  
     \\
   T^{x}(M). T^{x}(T^{r}(V)) &   & (D5)  & & (D6)   & T^{xrx^{-1}}(M). T^{xrx\inv}(T^{x}( V))   \\
   M . T^{x}(T^{r}(V)) & & & & M . T^{xr}(V) & M . T^{xrx\inv}(T^{x}( V) )\\
     };
  \path[-stealth]
    (m-1-1) edge node [left] {$T^{x}((T^{r}_{2})^{M,V})$} (m-3-1)
    (m-1-3) edge [dashed] node [right] {$(T_{2}^{xr})^{M, V}$} (m-2-3)
    (m-1-1) edge node [above] {$(T_{2}^{x , r})^{M, V}$}(m-1-3)
  (m-3-1) edge node [left] {$T^{x}(\mu_{M}^{r}.1)$} (m-5-1)
  (m-1-3) edge node [right] [above]{$(T_{2}^{xrx\inv, x})^{-1}_{M. V}$} (m-1-6)
         (m-1-6) edge node [left] {$T^{xrx^{-1}}((T_{2}^{x})^{M, V})$} (m-3-6)
          (m-3-1) edge [dashed] node [above] {$(T_{2}^{x})^{T^{r}(M)\;T^{r}(V)}$}(m-3-3)
           (m-5-1) edge  node [left] {$(T^{x}_{2})^{M, T^{r}(V)}$}(m-6-1)
           (m-6-1) edge  node [left] {$\mu^{x}_{M}. 1$}(m-7-1)
            (m-7-1) edge  node [below] {$1. (T_{2}^{x, r})_{V}$}(m-7-5)
            (m-7-5) edge  node [below] {$1 . (T_{2}^{xrx^{-1}, x})^{-1}_{V}$} (m-7-6)
             (m-3-3) edge [dashed] node [right] {$T^{x}(\mu_{M}^{r}). 1$} (m-6-1)
    (m-2-3) edge [dashed]  node [left] {$(T_{2}^{x, r})^{-1}_{M}. (T_{2}^{x, r})^{-1}_{V}$} (m-3-3)
      (m-2-3) edge [dashed]  node [right] {$\mu^{xr}_{M}. 1$} (m-7-5)
      (m-2-3) edge [dashed]  node [right] {$(T_{2}^{xrx\inv, x})^{-1}_{M}. (T_{2}^{xrx\inv, x})^{-1}_{V}$} (m-4-4)
      (m-3-6) edge [dashed]  node [above] {$(T^{xrx^{-1}}_{2})^{T^{x}(M), T^{x}(V)}$} (m-4-4)
(m-4-4) edge [dashed] node [left] {$T^{xrx^{-1}}(\mu^{x}_{M}). 1$} (m-6-6)   
(m-5-6) edge node [left] {$(T_{2}^{xrx^{-1}})^{M, T^{x}(V)}$} (m-6-6)
(m-6-6) edge node [left] {$\mu_{M}^{xrx^{-1}}. 1$} (m-7-6)
(m-3-6) edge node [left] {$T^{xrx\inv}(\mu_{M}^{x}. 1)$} (m-5-6)
         ;
\end{tikzpicture}
}
\end{center}
The pentagon $(D1)$ from the upper left corner is commutative by the fact that $T_{2}^{x,r}$ is a natural isomorphism of monoidal functors. On the other hand the pentagon $(D2)$ from the upper right corner is commutative by the fact that $T_{2}^{xrx\inv,x}$ is a natural isomorphism of monoidal functors. Note that the commutativity of the diagram $(D3)$ follows by the naturality of the transformation $T_{2}^{x}$, i.e Equation \eqref{tensor-rho}. A similar argument applies for diagram $(D4)$. Commutativity of the diagrams $(D5)$ and $(D6)$ follows by Equation \eqref{deltau}. 
\epf\br\label{conjstr}
Note that by the first statement of the previous lemma the monoidal functor $c_{H, x}:\cc^{H}\ra \cc^{\;^{x}H}$ also becomes a $\cc^{G}$-module functor. Indeed, for any $M\in \cc^{G}$ and any $V \in \cc^{H}$ one has that
\beqn
c_{H, x}(M \bxt V)=c_{H, x}(\re^{G}_{H}(M)\ot V)\simeq c_{H, x}(\re^{G}_{H}(M))\ot c_{H, x}(V)=\re^{G}_{\;^{x}_{H}}c_{G, x}(M)\ot c_{H, x}(V)
\eeqn
On the other hand since $c_{G, x}\simeq \id_{\cc^{G}}$ it follows that $c_{H, x}(M \bxt V)\simeq M \bxt c_{H, x}(V)$, which shows that $c_{H, x}$ is a $\cc^{G}$-module functor. Its module structure it is given by
\beq
T^{x}(M\ot V)\xra{T^{x}_{2}}T^{x}(M)\ot T^{x}(V)\xra{\mu_{M}^{x}\ot 1}M \ot T^{x}(V).
\eeq
\er
\ncm{\wta}{With the above notations one has}
\subsection{Proof of Theorem \ref{macky}}
We are now ready to give a proof for Theorem \ref{macky}. 
\bpf  Fix a set $\mathcal{D}$ of representative elements  for the double cosets $K \backslash H \slash L$. For any $x \in \mtc D$ consider an arbitrary set $\mtc{R}_{x}$ of representative elements for the  left cosets of $yL$ of $KxL/L$. Since  $H=\sqcup_{x\in \mtc D}KxL$ it follows that 
$\mtc R:=\sqcup_{x\in \mtc D}\mtc{R}_{x}$
is a complete set of representative elements  for the left cosets $H/L$. Suppose that $(V, \mu_{V}) \in \cc^{L}$. Then the induction functor $\ind^{L}_{H}$ associated to $\mtc R$ can be written as
\beqn
 \ind^{H}_{L}(V)\simeq (\bigoplus_{r \in \mathcal{R}}T^{r}(V), 
 \;\mu_{ \ind^{G}_{H}(V)}) = (\bigoplus_{x \in \mtc D}\bigoplus_{a \in \mtc R_{x}}T^{a}(V), \mu_{ \ind^{G}_{H}(V)})
 \eeqn
 where by Equation \eqref{muind} the equivariant structure is  given on components by
{
 \beqn
 \mu^{g, a}_{ \ind^{H}_{L}(M)}:T^{g}(T^{a}(V))\xra{(T_{2}^{g,a})_{V}} T^{ga}(V)=T^{yh}(V)\xra{(T_{2}^{y,l})^{-1}_{V}} T^{y}T^{l}(V)\xra{T^{y}(\mu_{Y}^{l})} T^{y}(V)
 \eeqn
 }
 if $ga=yl$ with $y \in \mathcal{R}$ and $l \in L$.

\noindent
For any $x \in \mtc D$ let 
 \beq 
F_x(V):=  \bigoplus_{a \in \mtc R_{x}}T^a(V).
 \eeq
It can be easily verified that the above induced equivariant structure $\mu_{\ind^{G}_{H}(V)}$ of $\ind^{G}_{H}(V)$ sends the component $F_x(V)$ to itself. Indeed, note that for any $a \in \mtc R_{x}$ if $m \in K$ and $ma=bl$ with $l \in L$ then $b=mal\inv \in KaH=KxH$. It follows that  $(F_x(V), \nu|_{K}) \in \cc^{K}$ and then one can define a functor $F_{x} :\cc^{L}\ra \cc^{K}$. Moreover it is easy to see that $
\res_K^H\ind^H_L\simeq\bigoplus_{x \in \mathcal D}F_x
$ as $k$-linear functors.

\noindent
Define also the functor  $G_{x}:\cc^{\;^{x}L}\ra \cc^{K}$ given by $G_{x}=\ind_{K \cap \;^xL}^K\circ \res^{\;^xL}_{\;^xL\cap K}$. Then in order to finish the proof it is enough to show that  for all $x \in D$ one has
\beq\label{pieces}
F_{x}\simeq G_{x}\circ c_{L, x}.
 \eeq
\noindent Note that there is a bijection between the following sets of left cosets
$
K/K\cap \;^xL \xra{\phi} KxL/L
$
 given by $a(K\cap \;^xL)\mapsto axL$ whose inverse is given by $aL\mapsto ax^{-1}(K\cap \;^xL)$. This enables us to write
$
G_{x}(P)\simeq \oplus_{a \in \mtc R_{x}}T^{ax^{-1}}(P),
$
for any $P \in \cc^{\;^{x}L}$. Under this isomorphism  the $K$-equivariant structure of $G_{x}(P)$ becomes
\beqarn
&& \mu_{G_{x}(P)}^{m}:T^{m}(T^{ax^{-1}}(P)) \xra{(T_{2}^{m , ax^{-1}})_{P}}T^{max^{-1}}(P) =T^{bx^{-1}{m'}}(P) \xra{}  \\ &&  \xra{(T_{2}^{bx^{-1} , m'})^{-1}_{P}}T^{bx^{-1}}(T^{{m'}}(P))\xra{T^{bx^{-1}}(\mu_{P}^{m'})} T^{bx^{-1}}(P)
\eeqarn
where $b \in \mtc R_{x}$ is chosen such that $max^{-1}=bx^{-1}{m'}$ with ${m'} \in \;^{x}L\cap K$. Thus for any $V \in \cc^{L}$ the $K$-equivariant structure of $G_{x}c_{L, x}(V)$ is given by
{\small
\beqarn
&& \mu_{G_{x}c_{H, x}(V)}^{m}:T^{m}(T^{ax^{-1}}(T^x(V))) \xra{(T_{2}^{m , ax^{-1}})_{T^x(V)}}T^{max^{-1}}(T^x(V)) =T^{bx^{-1}{m'}}(T^x(V)) \xra{}  \\ &&  \xra{(T_{2}^{bx^{-1} , m'})^{-1}_{T^x(V)}}T^{bx^{-1}}(T^{{m'}}(T^x(V)))\xra{T^{bx^{-1}}(T_{2}^{m',x})_{V}} 
T^{bx^{-1}}(T^{m'x}(V))\xra{} \\& & \xra{T^{bx^{-1}}((T_{2}^{x, x\inv m' x})^{-1}_{V})} T^{bx^{-1}}(T^{x}(T^{x^{-1}m'x}(V))) \xra{T^{bx^{-1}}(T^{x}(\mu_{V}^{x^{-1}m' x}))} T^{bx^{-1}}(T^{x}(V)).
\eeqarn
}
\noindent Define the natural transformation $N_{x}:F_{x}\ra G_{x}\circ c_{L, x}$ by \beqn (N_{x})_{V}: F_{x}(V) \xra{\bigoplus_{a \in \mtc R_{x}}(T^{ax^{-1}, x}_{2})^{-1}_{V} }G_{x}(c_{L, x}(V))\eeqn 
\noindent
In order to show that $N_{x}$ is well defined  it is enough to check that it is compatible with the two equivariant structures of the objects $F_{x}(V)$ and $G_{x}c_{H, x}(V)$. Thus one has to verify that the following diagram is commutative: 
\begin{center}
{\tiny
\begin{tikzpicture}
  \matrix (m) [matrix of math nodes, row sep=6em, column sep=6em,minimum width=5em]
  { T^{m}(F_{x}(V)) & F_{x}(V) \\
    T^{m}(G_{x} c_{L,x}(V)) & G_{x}c_{L,x}(V) \\
    };
  \path[-stealth]
  (m-1-1) edge node [above] {$\mu^{m}_{F_{x}(V)}$} (m-1-2)
    (m-2-1) edge node [above] {$\mu^{m}_{G_{x}c_{L,x}(V)}$} (m-2-2)
         
               (m-1-1) edge node [left] {$T^{m}((N_{x})_{V})$} (m-2-1)
    (m-1-2) edge node [right] {$(N_{x})_{V}$} (m-2-2)
   ;
\end{tikzpicture}
}
   \end{center}
\noindent On the components the above diagram is equivalent to the following. Suppose that $m \in K$ and $a, b \in \mtc R_{x}$ with $max^{-1}=bx^{-1}m'$ for some $m' \in K\cap \;^{x}L$. Then $ma=bx^{-1}m'x=bl$ with $l=(x\inv m'x) \in L$.
\begin{center}
{\tiny
\begin{tikzpicture}
  \matrix (m) [matrix of math nodes, row sep=6em,column sep=0.7em,minimum width=4em]
  {T^{m}(T^{ax^{-1}}(T^{x}(V))) & T^{max^{-1}}(T^{x}(V))=T^{bx^{-1}m'}(T^{x}(V)) & & T^{bx^{-1}}(T^{m'}(T^{x}(V)))  \\
  T^{m}(T^{a}(V))  &  &  & & \\ 
   T^{ma}(V)=T^{bl}(V) & & & T^{bx^{-1}}(T^{m'x}(V))=T^{b}(T^{xl}(V))\\
  T^{b}(T^{l}(V)) &  & & T^{bx^{-1}}(T^{x}(T^{l}(V)))\\
  T^{b}(V) &  & &T^{bx\inv}(T^{x}(V)) \\
  \\
};
  \path[-stealth] 
  (m-1-1) edge node [left] {$T^{m}[(T_{2}^{ax^{-1},x})_{V}]$} (m-2-1) 
  (m-3-4) edge node [right] {$T^{bx\inv}(T_{2}^{x,l})_{V}$} (m-4-4)
   (m-4-4) edge node [right] {$T^{b}(T^{x}(\mu^{l}_{V}))$} (m-5-4)
      (m-5-4) edge node [above] {$(T_{2}^{bx^{-1}, x})^{-1}_{V}$} (m-5-1)
  (m-2-1) edge node [left] {$(T_{2})^{m,a}_{V}$} (m-3-1)
  (m-1-4) edge  node [left] {$T^{bx\inv}(T^{m', x}_{2})_{V}$} (m-3-4)
   (m-1-2) edge [dashed] node [right]{$(T^{max\inv,x}_{{2}})_{V}=(T^{bx\inv m',x}_{{2}})_{V}$} (m-3-1) 
   (m-4-1) edge node [left]{$T^{b}(\mu^{l}_{V})$} (m-5-1)

    (m-1-1) edge node [above]{$(T^{m,ax^{-1}}_{{2}})_{T^{x}(V)}$} (m-1-2) 
     (m-1-2) edge node [above]{$(T^{bx^{-1}, m'}_{{2}})_{T^{x}(V)}$} (m-1-4) 
      (m-3-1) edge [dashed]  node [above]{$(T_{2}^{bx\inv, xl})_{V}$} (m-3-4) 
           (m-3-1) edge  node [left]{$(T_{2}^{b, l})^{-1}_{V}$} (m-4-1) 
           (m-4-1) edge [dashed]  node [above]{$(T_{2}^{bx\inv,x})^{-1}_{T^{l}(V)}$} (m-4-4) 
    ;
\end{tikzpicture}
}
\end{center}
\noindent
The bottom rectangle is commutative by the naturally of $T^{bx\inv, x}_{2}$ with respect to the morphisms, Equation \eqref{nat-ro}. The above rectangle is commutative due to Equation \eqref{ro-2}, the associativity of the action. The upper left diagram is commutative by the same reason. The upper right trapeze is commutative by associativity of the action,  Equation \eqref{ro-2}. Clearly  $(N_{x})_{V}$ is an isomorphism in $\cc$ and therefore it is an isomorphism in $\cc^{K}$.

\noindent 2) Suppose now that $\cc$ is a $k$-linear monoidal category  and the action of $G$ is by monoidal autoequivalences. It remains to show that  the  isomorphism from the statement is an isomorphism of $\cc^{G}$-module functors.  For this, it is enough to show that that  the above isomorphism $N_{x}$ of Equation \eqref{pieces} is an isomorphism of $\cc^G$-module functors. This is equivalent to the commutativity of the following diagram:
{
\begin{center}
{\tiny
\begin{tikzpicture}
  \matrix (m) [matrix of math nodes, row sep=9em,column sep=6em,minimum width=5em]
  { 
G_{x}c_{H,x}(M\boxtimes V) &
M\boxtimes G_{x}c_{H, x}(V) \\
   F_{x}(M\boxtimes V) &
M\boxtimes F_{x}(V).\\
    };
  \path[-stealth]
  (m-1-1) edge node [above] {$(G_{x}c_{H,x})^{M, V}_{2}$} (m-1-2)
    (m-2-1) edge node [above] {$(F_{x})^{M, V}_{2}$} (m-2-2)
               (m-1-1) edge node [left] {$(N_{x})_{M\boxtimes V}$} (m-2-1)
    (m-1-2) edge node [right] {$1_{M}\boxtimes (N_{x})_{ V}$} (m-2-2)
   ;
\end{tikzpicture}
}
   \end{center}
}\noindent 
Note that the $\cc^{G}$-module structure of $F_{x}$ is the one induced from $\ind^{H}_{L}$.  Thus one has that $(F_{x})_{2}^{M, V}:F_{x}(M \otimes V)\ra M \otimes F_{x}(V)$ is defined on the components as
\beqn
(F_{x})_{2}^{a, M, V}:T^a(M \ot V) \xra{{(T^a)_{2}}^{M, V}} T^a(M) \ot T^{a}(V) \xra{\mu^{a}_{M}\ot 1}M \ot T^{a}(V). 
\eeqn
\noindent
On the other hand note that the $\cc^{G}$-module functor structure of  $G_{x}\circ c_{H, x}$ on the components is given by:
{
\beqarn
&& (G_{x}c_{H, x})^{a, M, V}_{2}:T^{ax^{-1}}(T^{x}(M \ot V))\xra{T^{ax^{-1}}(T^{x})_{2}^{M, V}} T^{ax^{-1}}(T^{x}(M) \ot T^x (V))\xra{}
\\ && \xra{T^{ax^{-1}}(\mu_{M}^{x}\ot 1) }T^{ax^{-1}}(M \ot T^{x}(V)) \xra{(T^{ax^{-1}})_{2}^{M, T^x(V)}} T^{ax^{-1}}(M)\ot T^{ax^{-1}}(T^{x}(V)) \xra{} \\ && \xra{\mu_{M}^{ax^{-1}}\ot 1} M \ot T^{ax^{-1}}(T^x(V)).\eeqarn
}
\noindent
Thus for any $a \in \mtc R_{x}$, on the components, the above diagram is equivalent to the following:
\begin{center}
{\Tiny
\begin{tikzpicture}
  \matrix (m) [matrix of math nodes, row sep=4em,column sep= 0.3 em,minimum width=3em]
  { T^{ax^{-1}}(T^{x}(M\ot V)) &   & T^{a}(M\ot V)\\
   T^{ax^{-1}}(T^{x}(M)\ot T^{x}(V)) & T^{ax^{-1}}(T^{x}(M))\ot T^{ax^{-1}}(T^{x}(V)) & T^{a}(M)\ot T^{a}(V)\\
     T^{ax^{-1}}(M\ot T^{x}(V)) & T^{a}(M)\ot T^{ax^{-1}}(T^{x}(V)) & \\
  T^{ax^{-1}}(M)\ot T^{ax^{-1}}(T^{x}(V)) & & \\
 M \otimes  T^{ax^{-1}}(T^{x}(V)) & &  M \ot T^{a}(V)\\
 };
  \path[-stealth]
  (m-1-1) edge node [above] {$(T^{a, x}_{2})_{M\ot V}$} (m-1-3)
         (m-4-1) edge node [left] {$\mu_{M}^{ax^{-1}}\ot \id $} (m-5-1)
               (m-1-1) edge node [left] {$T^{ax^{-1}}((T^{x}_{2})_{M, V})$} (m-2-1)
                (m-2-1) edge node [left] {$T^{ax^{-1}}(\mu_{M}^{x}\ot 1)$} (m-3-1)
                 (m-2-1) edge [dashed]   node [above] {$$} (m-2-2)
                 (m-2-2) edge [dashed]   node [above] {$$} (m-2-3)
                  (m-2-2) edge [dashed]   node [above] {$$} (m-3-2)
                  (m-2-2) edge [dashed]   node [right] {$ $} (m-4-1)
                  (m-3-2) edge [dashed]   node [above] {$$} (m-5-1)
                  (m-3-1) edge node [left] {$(T^{ax^{-1}, x})_{2}^{M, T^{ax^{-1}}(V)}\ot \id$} (m-4-1)
                    (m-2-3) edge node [right] {$\mu_{M}^{a}\ot \id$} (m-5-3)
    (m-1-3) edge node [right] {$(T^{a}_{2})^{{M, V}}$} (m-2-3)
     (m-5-1) edge  node [above] {$\id \ot (T_{2})^{ax^{-1},x}_{V}$} (m-5-3)
    ;
\end{tikzpicture}
}
\end{center}
\noindent
Note that the upper rectangle is commutative since $T^{ax,x^{-1}}_{2}$ is a natural transformation of monoidal functors. The part below is commutative by the natural properties of the tensor bifunctor of the monoidal  category $\cc$ and the compatibility between $\mu_{M}^{ax^{-1}}, \mu_{M}^{x^{-1}}$ and $\mu_{M}^{a}$.
\epf
\br Note that $\Rep(G)$ can be regarded as the equivariantization $\Vc^{G}$ of the trivial action of $G$ on $\cc=\Vc$, \cite{ENO}. In this case the previous theorem recovers the usual Mackey decomposition for representations of finite groups.
\er
\bl\label{compind} Suppose that a finite group $G$ acts $k$-linearly on $\cc$. Let $K\leq H \leq L$ be a tower of subgroups of $G$. 
1) One has that
 \beqn
 \mtc R^{L}_{K} =\mtc R^{H}_{K} \mtc R^{L}_{H}  
\eeqn

\noindent
2) There is a natural transformation \beqn
{\bf I}^{H}_{K,L}:\mtc I^{L}_{K}\ra \mtc I^{L}_{H}\mtc I^{H}_{K}
\eeqn which is an isomorphism of $k$-linear functor.

\noindent
3) Moreover if the action of $G$ on $\cc$ is by monoidal equivalences then ${\bf I}^{H}_{K,L}$ is an isomorphism of $\cc^{L}$-module functors.
\el 
\bpf The natural transformation ${\bf I}^{H}_{K,L}$ is defined as follows. Fix $\mtc R$ and $\mtc S$ sets of representative elements  for the left cosets of $K$ inside $H$ and of $H$ inside $L$ respectively. Then $\mtc R\mtc S:=\{rs\;|\; \in R,\; s\in S\}$ is a set  of representative for the left cosets of $K$ inside $L$. Then for any $M \in \cc^{K}$ one has
$\mtc I^{K}_{L}(M)=\bigoplus_{r \in \mtc R\;s \in \mtc S}T^{rs}(M)$ while $\mtc I^{K}_{H}(M)=\bigoplus_{r \in \mtc R}T^{r}(M)
$ and $\mtc I^{H}_{L}(P)=\bigoplus_{s \in \mtc S}T^{s}(P)
$ for any $P \in \cc^{H}$. Then one can define on the components
\beqn
({\bf I}^{H}_{K,L})_{M}=\bigoplus_{r \in \mtc R\;s \in \mtc S}(T^{r,s}_{2})_{M}.
\eeqn
Similarly to Lemma \ref{commir} one can check that $({\bf I}^{H}_{K,L})_{_{M}}$ is an isomorphism in $\cc^{L}$. Thus one has to verify that the following diagram {\small
\begin{center}
{\tiny
\begin{tikzpicture}
  \matrix (m) [matrix of math nodes, row sep=9em,column sep=8em,minimum width=5em]
  { T^{l}(\ind^{L}_{K}( M) ) & T^{l}(\ind^{L}_{H}\ind^{H}_{K}( M)) \\
 \ind^{L}_{K}( M) & \ind^{L}_{H}\ind^{H}_{K}( M)\\
    };
  \path[-stealth]
  (m-1-1) edge node [above] {$T^{l}(({\bf I}^{H}_{K,L})_{M})$} (m-1-2)
    (m-2-1) edge node [above] {$({\bf I}^{H}_{K,L})_{M}$} (m-2-2)
         (m-1-1) edge node [left] {$\mu^{l}_{\ind^{L}_{K}( M)}$} (m-2-1)
    (m-1-2) edge node [right] {$\mu^{l}_{\ind^{L}_{H}\ind^{H}_{K}( M)}$} (m-2-2)
   ;
\end{tikzpicture}
}
   \end{center}
}\noindent
is commutative for any $l\in L$.
\noindent Indeed suppose that $lsr=s'r'k'$ for some $k \in K$ and $s'\in \mtc S$ and $r' \in \mtc R$. Morever suppose that $ls=s''h$ and $hr=r''k''$ for some $r'' \in \mtc R$, $s''\in \mtc S$, $h \in H$ and $k''\in K$. Since $lsr=s'r'k'=s''r''k''$ it follows that $s'=s''$, $r'=r''$ and $k'=k''$. Thus, on the components the above diagram is equivalent to the following diagram:
\begin{center}
{\Tiny
\begin{tikzpicture}
  \matrix (m) [matrix of math nodes, row sep=7em,column sep= 5.9 em,minimum width=3em]
  { T^{l}(T^{sr}(M) ) &    T^{l}( T^{s}( T^{r}(M)))\\
   T^{lsr}(M) & T^{ls}(T^{r}(M))\\
     T^{s'r'}(T^{k'}M) & T^{s'}( T^{h}( T^{r}(M))) \\
  T^{s'r'}(M)& T^{s'}( T^{hr}(M))  \\
 &  T^{s'}( T^{r'}( T^{k'}(M)))\\
 & T^{s'}( T^{r'}(M))\\
 };
  \path[-stealth]
         (m-1-1) edge node [above] {$T^{l}((T^{s, r}_{2})^{-1}_{M})$} (m-1-2)
         (m-1-2) edge node [right] {$(T^{a}_{2})^{{M, V}}$} (m-2-2)
         (m-4-2) edge node [right] {$T^{s'}((T_{2}^{r', k'})^{-1}_{M} )$} (m-5-2)
         (m-1-1) edge node [left] {$(T^{l,sr}_{2})_{M}$} (m-2-1)
         (m-2-1) edge node [left] {$(T^{s'r', k'}_{2})_{M}$} (m-3-1)
         (m-2-1) edge [dashed]   node [above] {$(T^{ls, r}_{2})^{-1}_{M}$} (m-2-2)
         (m-2-1) edge [dashed]   node [above] {$(T^{s',hr}_{2})_{M}^{-1}$} (m-4-2)
         (m-2-2) edge node [right] {$(T_{2}^{s', h})^{-1}_{T^{r}(M)}$} (m-3-2)
         (m-5-2) edge node [right] {$T^{s'}(T^{r'}(\mu_{M}^{k'}))$} (m-6-2)
         (m-3-2) edge  node [right] {$T^{s'}((T_{2}^{h,r})_{M})$} (m-4-2)
         (m-3-1) edge node [left] {$ T^{s'r'}(\mu^{k'}_{M})$} (m-4-1)
          (m-3-1) edge [dashed]   node [above]  {$(T^{s',r'}_{2})_{T^{k}(M)}^{-1}$} (m-5-2)
         (m-4-1) edge [dashed]   node [below] {$ (T_{2}^{s',r'})^{-1}_{M}$} (m-6-2)
    ;
\end{tikzpicture}
}
\end{center}
\noindent
2) Suppose now that $\cc$ is a $k$-linear monoidal category  and the action of $G$ is by monoidal autoequivalences. Let $M \in \cc^{L}$ and $r \in \mtc R$, $s \in \mtc S$. The fact that $\mtc I^{K}_{H}$ is natural transformation of $\cc^{L}$-module functors follows from the commutativity of the diagram below.
\begin{center}
{\tiny
\begin{tikzpicture}
  \matrix (m) [matrix of math nodes, row sep=3.25 em,column sep=3 em,minimum width=0.1em]
  { 
     T^{rs}(M\ot V)) &     & T^{r}(T^{s}(M\ot V))\\
     &  & T^{r}(T^{s}(M)\ot T^{s}( V))\\
    T^{rs} (M)\ot  T^{rs}(V) &  T^{r}(T^{s}(M)) \ot T^{r}(T^{s}(V))  &  \\
     M\ot T^{rs}(V) & T^{r}(M) \ot T^{r}(T^{s}(V)) & T^{r}(M\ot T^{s}(V))\\
     };
  \path[-stealth]
    (m-1-1) edge node [left] {$(T^{rs}_{2})^{M, V}$} (m-3-1)
    (m-3-1) edge [dashed] node [above] {$(T_{2}^{r, s})^{-1}_{M}\ot (T_{2}^{r, s})^{-1}_{V} $} (m-3-2)
    (m-1-1) edge node [above] {$(T^{r,s}_{2})_{M \ot V}$}(m-1-3)
  (m-1-3) edge node [right] {$T^{r}((T^{s}_{2})^{M,V})$} (m-2-3)
  (m-2-3) edge node [right]{${T^{r}(\mu^{s}_{M}\ot 1)}$} (m-4-3)
         (m-4-3) edge node [above] {$(T^{r}_{2})^{M, T^{s}(V)}$} (m-4-2)
          (m-4-2) edge  node [above] {$\mu_{M}^{r}\ot (T_{2}^{r,s})_{V}$}(m-4-1)
           (m-3-1) edge node [left] {$\mu_{M}^{rs}\ot 1$}(m-4-1)
            (m-3-2) edge [dashed] node [right] {${T^{r}(\mu^{s}_{M})\ot 1}$}(m-4-2)
             (m-2-3) edge [dashed] node [left] {$(T^{r}_{2})^{T^{s}(M), T^{s}(V)}$}(m-3-2)
         ;
\end{tikzpicture}
}
\end{center}
The upper pentagon is commutative by Equation \eqref{tensor-rho}. The bottom left square is commutative by Equation \eqref{deltau}.  The bottom right square is commutative by the naturality of the transformation $T^{r}_{2}$ with respect to the first argument.
\epf
\bp\lb{atp}Suppose that $\cc, \cd$ and $ \ce$ are rigid monoidal categories and $F_1:\cc \ra \cd$ and $F_2:\cc \ra \ce$ are two monoidal functors with left adjoint functors $I_1:\cd\ra \cc$ and respectively $I_2:\ce\ra \cc$.  Then for any objects $M \in \co(\cd)$ and $N \in \co(\ce)$ one has the canonical isomorphism in $\cc$
\beq\lb{ti}
I_1(M)\ot I_2(N)\simeq I_1(F_{1}(I_2(N))\ot M) \simeq I_2(F_2(I_1(M))\ot N).
\eeq

\ep
\bpf
It can be shown by a straightforward computation that 
\beq\lb{ti}
\Hom_{\cc}(I_1(M)\ot I_2(N), P)\simeq \Hom_{\cc}(I_2(F_2(I_1(M))\ot N)
, P)
\eeq
for any object $P \in \cc$. Indeed, 
\begin{eqnarray*}
  \Hom_{\C}(I_1(M)\ot I_2(N), P)&=& \Hom_{\C}(I_{1}(M),\;P \ot I_{2}(N)^*)= \\
  =\Hom_{\mtc{D}}(M,\;F_{1}(I_{2}(N)^*\ot P)) &=&  \Hom_{\mtc{D}}(M, F_{1}(I_{2}(N))^* \ot F_{1}(P))=\\
  =\Hom_{\mtc{D}}(F_{1}(I_{2}(N)) \ot M, F(P)) &=&  \Hom_{\mtc{C}}(I_1(F_{1}(I_2(N))\ot M), P)
\end{eqnarray*}

Then Yoneda's lemma  implies the conclusion.\epf

In particular for $\ce=\cc$ and $F_{2}=I_{2}=\id_{\cc}$  one obtains that 
\beq\lb{ti2}
I_{1}(M)\ot V\simeq I_{1}(M\ot F_{1}(V))
\eeq for any objects $M\in \cd$ and $V\in \cc$.
\subsection{Proof of Theorem \ref{cmain}.} 1) First we will define the functors $\mtc R^{L}_{H}$, $\mtc I^{L}_{H}$ and $c_{H,a}$.

\noindent 
The functors $\mathcal R_{L}^{H}$ correspond to the restriction functors defined in section \ref{indc}. We fix an arbitrary set of representative elements for the left cosets $Hx$ of any group inclusion $H\leq L$. We consider the corresponding induction functor  $\ind^{L}_{H}$ as in Equation \eqref{indx}. Moreover, the conjugation functors $c_{H,a}$ are defined as in Lemma \ref{conj}.

\noindent Secondly we will define the natural transformation from the definition of the categorical Mackey functor.
Note that the natural transformation ${\bf R}^{H}_{K, L}$ can be taken the identity by the first item of Lemma \ref{compind}. Also one can take  ${\bf CR }^{K}_{H, a}$ as identity functors by Equation \eqref{cr} of Lemma  \ref{commir}.  By Proposition \ref{nat-c} one can define $({\bf C}^{H}_{a,b})_{M}:=(T^{a,b}_{2})^{-1}_{M}$ for any object $M \in \cc^{\;^{ab}H}$.  Moreover the natural transformation
$
{\bf I}^{H}_{K,L}:\mtc I^{L}_{K}\ra \mtc I^{L}_{H}\mtc I^{H}_{K}
$ is defined by Lemma \ref{compind}. The definition of the natural transformations ${\bf CI}^{K}_{H,a}$ is given in Lemma \ref{commir}.

\noindent Now one has to verify all the compatibility conditions from the definition of a categorical Mackey functor. The Mackey  decomposition from Equation \eqref{mackeyss} is proven in Theorem \ref{macky}. Clearly the identities \eqref{idS}, \eqref{Sid} and \eqref{unul} hold. Moreover the diagrams $(\mathrm{R}), (\mathrm{RCC}), (\mathrm{RRC})$ are commutative since all  the natural transformations are the identity functors. Diagram $(C)$ is commutative by Equation \eqref{ro-2}. The commutativity of the diagrams $(\mathrm{ICC}), (\mathrm{IIC}), (\mathrm{I})$ follow by straightforward computation taking care of the representative elements for the left cosets of any inclusion of subgroups of $G$.

\noindent 2) Suppose now that the action of $G$ is by monoidal autoequivalences on the $k$-linear monoidal category  $\cc$. Note that we have already noticed that the restriction functors $\re^{H}_{L}$ are monoidal functors. By the second part of Lemma \ref{conj} one has that $c_{H, a}$ are also monoidal functors and the natural transformations ${\bf C}^{H}_{a,b}$ are morphisms of monoidal functors. 
The natural transformations ${\bf R}^{H}_{K, L}$ and ${\bf CR }^{K}_{H, a}$ are automatically morphisms of module functors since they are the identity functors. Note that the condition $\mathrm{(CG6)}$ follows from Proposition \ref{atp}.  All the other additional compatibility structures required for a categorical Green functor follow from the above lemmata.\qed

 \bt\label{kpass} Let $G$ be a group and $\mathcal{M}$ be a categorical $G$-Mackey functor over the monoidal category $\mtc S$. With the above notations one has the following:

\noindent
1)   Any categorical Mackey functor $\ct$ gives a usual Mackey functor $L \mapsto K_{i}(\ct(L))$ by taking the $K$-theory over $K_{0}(\cs)$.
  
\noindent
2) Moreover, if $\ct$ is a categorical Green functor then $K_{0}$ gives a Green functor $L \mapsto K_{0}(\ct(L))$ over $K_{0}(\cs)$.
  \et
\bpf
Similarly to \cite{kevin} one can use the following elementary facts about K-theory (see \cite{quillen}): If
$F_{1}$ and $F_{2}$ are isomorphic exact functors on an exact category, then they induce the
same map on $K$-theory;  and if $F_{1}$ and $F_{2}$ are exact functors on an exact category inducing
homomorphisms $f_{1}$ and $f_{2}$ on $K$-groups, then the functor $F_{1}\oplus F_{2}$ induces the
homomorphism $f_{1}+f_{2}$.
Now identities $\mathrm{(M1)}$ - $\mathrm{(M4)}$ from the definition of a Mackey functor follow from their functorial counterparts given in the definition of the categorical Mackey functor . Moreover, if $F$ is a monoidal functor then it induces an algebra morphism at the level of the Grothendieck rings. The adjunction properties from the definition of a Green functor follow from the Proposition \ref{atp}.
\epf
Applying Theorem \ref{kpass} to the Mackey functor fromTheorem \ref{cmain} we obtain Corollary \ref{main3}. We finish the paper with the following two examples previoulsy considered in the literature.
\bn{example}
Suppose that $R \subset S$ is a Galois extension of rings with Galois Group $G$. Then as in Example \ref{automequiv} the group $G$ acts on the category $S$-mod and $(S$-mod $)^{G}\simeq S\# \Z G$-mod. Since $S^{G}$ is Morita equivalent to $S\# \Z G$ and the $K$-theory is preserved by Morita equivalence it follows by Corollary \ref{main3} that $H \mapsto K_{i}(S^{H})$ is a Mackey functor. Thus our results extends the results obtained in \cite{kevin} for Galois extensions of commutative rings.
\end{example}
\bn{example}
Suppose that we have a cocentral extension of semisimple Hopf algebras 
\beq\lb{cocen}
k \ra B\xrightarrow{i} H\xrightarrow{\pi} kF\ra k.
\eeq
\noindent
Recall that this means the above sequence is exact, (see \cite{schss}), and  that $kF^*\subset \mtc{Z}(H^*)$ via $\pi^*$. 
Following \cite[Proposition 3.5]{natalecoc} it follows that $F$ acts on the fusion category $\Rep(B)$ and $\Rep(H)=\Rep(B)^F$.  Recall from \cite{natalecoc} that for all $x \in F$ the action is given by $T^x(M)=M$ as vector spaces with the action of $B$ on $T^x(M)$ given by
$
b.\;^xm=(x^{-1}.b)m.
$ Note that there is a typo in defining this action in \cite[Section 3.2]{natalecoc}. With this categorical group action the Green functors from Corollary \ref{main3} coincide to those described in \cite[Theorem 5.8]{jlms}. 
\end{example}
\bibliographystyle{amsplain}
\bibliography{gfr}
\ed
\section{Coherent group actions on graded  fusion categories}\lb{hmg} In this chapter we apply the results of the previous section to the case of fusion categories.  Let $\cc$ be a graded fusion category by a finite group $G$. Recall that this means that
$
\cc=\oplus_{g \in G}\cc_{g}
$
as abelian categories, and the tensor functor $\ot: \cc \times \cc\ra \cc $ sends $\cc_g\ot \cc_h$ into $\cc_{gh}$. For an object $V \in \cc$ define by $V_g$ the homogenous component of $V$ of degree $g$ from the above grading.

\noindent
Suppose further that another finite group $F$ acts by group automorphisms on $G$. Suppose that $F$ also acts by tensor automorphisms on the category $\cc$ via the action $T:F\ra \underline{\mtr{Aut}}_{\ot}(\cc)$ given by $x \mto T^x  :\cc \ra \cc$.

\bn{defn} An action of a finite group $F$ on the $G$-graded fusion category $\cc$ is called {\it coherent} with respect to the action of $F$ on $G$
 if 
\beq
T^x(\cc_g)\subset \cc_{\;^xg}
\eeq
for all $x \in F$ and $g \in G$.
\end{defn}
\mdn 
\bp Suppose that a finite group $F$ acts by tensor automorphisms on a $G$-graded fusion category $\cc$ via the action $T:F\ra \underline{\mtr{Aut}}_{\ot}(\cc)$ given by $x \mto T^x  :\cc \ra \cc$. Suppose further that $\cc_1$ is stable under the action of $F$. Then there is an action of $F$ on $G$ by group automorphisms such that the action of $F$ on $\cc$ is coherent with respect to this action.  
\ep
\bpf 
Let $V, W\in \cc_g$ and suppose that $T^x(V)\in \cc_{g_1}$ and $T^x(W) \in \cc_{g_2}$. Then $T^x(V^*\ot W)\in \cc_1$ since $V^*\ot W \in \cc_1$. On the other hand $T^x(V^*\ot W)\simeq T^x(V^*)\ot T^x(W)\in \cc_{{g_1}^{-1}g_2}$ which implies that $g_1=g_2$. Denoting $\;^xg:=g_1$ then it is easy to check that this defines an action of $F$ on $G$ by group automorphisms.
\epf
\bp
Suppose that a finite group $F$ acts on the fusion category $\cc$ and let $\cc=\oplus_{g \in U_{\cc}}\;\cc_{g}$ be the universal grading of $\cc$. Then:
\bne
\item The group $F$ acts on the group $U_{\cc}$ by automorphisms.
\item The action of $F$ on $\cc$ is coherent with respect to the action of $F$ on the group $U_{\cc}$.
\ene
\ep
\bpf It is known that $\cc_{1}:=\cc_{ad}$, the adjoint subcategory of $\cc$ it is generated by $V\ot V^{*}$ with $V$ a simple object of $\cc$. By the previous proposition it is enough to show that $T^{x}$ fixes $\cc_{1}$ for any $x \in F$.  Note that $T^{x}(V\ot V^{*})\simeq T^{x}(V)\ot T^x(V)^{*} \in \cc_{ad}$, thus the proof is completed.
\epf
\subsection{Examples of coherent actions of groups and their equivariantized categories} In this subsection we give some examples of coherent group actions on fusion categories.
\bn{example}\lb{gbb}{\it Braided $G$-crossed categories.}
Recall \cite{tuv} that a braided $G$-crossed fusion category is a quadruple $(\cc, G, T, c)$, where $G$ is a finite group, $\cc$ is a tensor category with a (not necessarily faithful) $G$-grading $\cc=\oplus_{{g \in G}}\cc_{g}$
 and a tensor action $T : G \ra \underline{\mtr{Aut}}_{\ot}(\cc)$, $g \mapsto T^{g}$ satisfying $T^{g}(\cc_{h})\subseteq \cc_{ghg\inv}$. Moreover the crossed braiding $c$ is defined by $c(X, Y):X\ot Y \ra T^{g}(Y)\ot X$ for all $X\in \cc_{g}$ and $Y \in \cc$. The compatibility conditions that have to be satisfied by  this datum can be found for example in \cite{tuv, DGNO}. Note that in this case $F=G$ acts coherently on $\cc$ with respect to the action of $G$ on  itself given by conjugation.\end{example}
 \bn{example}\lb{twisted} 
Group actions on graded pointed fusion categories are always coherent. See \cite[Section 4]{naidu}. In particular,  it follows by \cite[Lemma 6.3]{naidu} that the representation category of a (twisted) quantum double of a finite group is the equivariantization of a coherent action.
\end{example}
 \subsubsection{Universal gradings for the category of representations of semisimple Hopf algebras}
Let $A$ be a semisimple Hopf algebra over an algebraically closed field $k$. It is well known that $\mtr{Rep}(A)$ is a fusion category. Moreover there is a maximal central Hopf subalgebra $K(A)$ of $A$ such that $\mtr{Rep}(A//K(A))$ coincides to $\mtr{Rep}(A)_{
_{ad}}$ the adjoint subcategory of $\rep(A)$, see \cite[Theorem 2.4]{GN}. Since $K(A)$ is commutative it follows that
$K(A)=kG^*$ where $G$ is the universal grading group of $\mtr{Rep}(A)$. The (universal) grading on $\rep(A)$ is given by  \bq\label{grdm}
\rep(A)_{g}=\{M \in \Irr(A)\;|\; p_{g}m=m \;\text{for all}\;m\in M\}.
\eq
Here $p_{g}\in k^{G}$ is the dual basis of group element basis $g \in G$.
\subsubsection{Cocentral extensions of semisimple Hopf algebras}\label{cocs}

Suppose that we have a cocentral extension of Hopf algebras 
\beq\lb{cocen}
k \ra B\xrightarrow{i} H\xrightarrow{\pi} kF\ra k.
\eeq

Recall that this means the above sequence is exact \cite{schss} and  $kF^*\subset \mtc{Z}(H^*)$ via $\pi^*$. 
On the other hand, using the reconstruction theorem from \cite{AD} it follows that 
$
H \simeq B\;^{\tau}\#_{\sg} \;kF
$
for some cocycle $\sg:B\ot B \ra kF$ and some dual cocycle $\tau:kF \ra B\ot B$ satisfying certain compatibility axioms.  In this case there is a weak action of $F$ on $B$ denoted by $f.b$ such that the multiplication and comultiplication on $H$ become
\bq\lb{mu}
(b\#_{\sg}f)(c\#_{\sg}g)=b(f.c)\sg(f,g)\#_{\sg}fg
\eq
and respectively
\bq\lb{com}
\D(b\#_{\sg}\bar{f})=(b_1\tau(\bar{f})_i\#_{\sg}\bar{f})\ot(b_2\tau(\bar{f})^i\#_{\sg}\bar{f}).
\eq
\bp\label{cocsp} Suppose that we have a cocentral extension of semisimple Hopf algebras as in Equation \eqref{cocen}.
Moreover, with the above notations suppose that $\cc:=\rep(B)$ and let
$\mtc{C}=\oplus_{g \in G}\mtc{C}_g
$
be the universal grading of $\rep(B)$ where $K(B)=kG^*$. Then $F$ acts coherently on $\cc$ with respect to a given action by group automorphisms of $F$ on $G$.
\ep
\bpf
First it will be shown that $F$ acts on $G$ by group automorphisms. In order to do this we show that $F$ acts on $K(B)=kG^*$ by Hopf automorphisms.
\mdn Note that if $b \in K(B)$ then $f.b \in K(B)$ where $``.``$ represents the weak action of $F$ on $B$ from above. Indeed $f.b \in \cz(B)$ since $F$ acts by algebra automorphisms on $B$.
On the other hand using formula $(A)$ from \blue{\cite[Section 2]{AD} } it follows that 
\beq
\Delta(f.b)=\tau(f)(f.b_1\ot f.b_2)\tau(f)^{-1}.
\eeq
Therefore if $b \in \cz(B)$ then 
$
\Delta(f.b)=f.b_1\ot f.b_2
$
which shows that $F. K(B)$ is a central Hopf subalgebra of $B$.  Thus $F.K(B)\subseteq K(B)$ and $F$ acts by Hopf algebra automorphisms on $K(B)$. This implies that there is an action of $F$ on $G$ such that the action of $F$ on $K(B)=kG^*$ is given by
$
x.p_g=p_{\;^xg}
$
for all $x \in F$ and $g \in G$.
\mdn Following \cite[Proposition 3.5]{natalecoc} it follows that $F$ acts on the fusion category $\Rep(B)$ and $\Rep(H)=\Rep(B)^F$.  Recall from \cite{natalecoc} that for all $x \in F$ the action is given by $T^x(M)=M$ as vector spaces with the action of $B$ on $T^x(M)$ given by
$
b.\;^xm=(x^{-1}.b)m.
$

\mdn In order to verify  that the above action of $F$ on $\rep(B)$ is coherent with respect to  this action one has to verify that if $M \in \co(\cc_g)$ then $T^x(M)\in \co( \cc_{\;^xg})$. Using Equation \eqref{grdm} it follows that for any $h \in G$ one has that 
${p_h.}\;^xm=(x^{-1}.p_h)m=p_{\;^{x^{-1}}h}\;m=\delta_{\;^{x^{-1}}h,\; g}\;m=\delta_{h, \;^{x}g}\;m$
which shows that indeed $T^x(M)\in \co(\cc_{\;^xg})$.
\epf
\noindent \bll{Put the result that for the universal grading any action is coherent. Make the above proposition also true for the nonsemisimple case.}

\sub{Simple objects for equivariantizations of coherent actions}In this subsection we investigate the simple objects of an equivariantization under a coherent action. 
\md

Suppose that $F$ acts coherently on a $G$-graded fusion category $\cc$ with respect to a given action of $F$ on $G$.
With the above notations note that the stabilizer $F_{g}$ of an element $g \in G$ acts by $k$-linear automorphisms on the abelian subcategory $\cc_{g}$ of $\cc$. In particular one obtains in this way an action of  $F$ on the fusion subcategory $\cc_{1}$ of $\cc$. Note that this action on $\cc_{1}$ is by tensor automorphisms.

\bl Suppose that $F$ acts coherently on a  $G$-graded fusion category $\cc$ with respect to an action by group automorphisms on $G$. Moreover suppose that $(V, {\mu_V^x}|_{\{x \in F\}})\in \cc^F$ is an equivariantized object with a canonical decomposition $V=\oplus_{g \in G}V_{g}$. Then for all $g \in G$ one has that $(V_g,  {\mu^x_{V_{g}}}|_{\{x \in F_g\}}) \in \cc^{F_g}$.
\el

\bpf If $V=\oplus_{g \in G}V_{g}$ then $T^{x}(V)=\oplus_{g \in G}T^{x}(V_{g})$.
Since $\mu_V^x:T^x(V)\ra V$ is an isomorphism in $\cc$ it follows that $\mu^x_V$ sends the component $T^{x}(V_g)$ of $T^{x}(V)$ to the component $V_{\;^{x}g}$ of $V$. Thus if $x \in F_g$ then $\mu^{x}_{V}|_{T^{x}(V_{g})}:T^{x}(V_g)\ra V_{g}$ is an isomorphism. Note that the compatibility conditions from Equation \eqref{deltau} for an equivariantized object $V \in \cc^{G}$ implies that  $V_{g} \in \cc^{F_{g}}$. \epf

Let $\Gm$ be a set of representative elements  for the orbits of the action of $F$ on the group $G$. For any $g \in G$ let $\co(g)$ be the orbit of $g$ under the action of $F$. Next proposition is a generalization of \cite[Proposition 2.7]{gnn}.

\bp\label{desgnn} Suppose that $F$ acts coherently on a  $G$-graded fusion category $\cc$ with respect to an action by group automorphisms on $G$. Then the set of simple objects of $\cc^F$ is parametrized by pairs $(g , M)$ where $g \in \Gm$ and $M \in {\cc_g}^{F_g}$  is a simple object. The simple object associated to the pair $(g, M)$ is given by the induced object $\Ind_{F_{g}}^{F}(M)$.
\ep 

\bpf  
Suppose that $(V, {\mu^{x}_{V}}_{\{x \in F\}})\in \cc^F$ is a simple object and let $V=\oplus_{g \in G}V_g$ with $V_{g }\in \cc_{g}$ be its decomposition viewed as an object of $\cc$. If $V_g \neq 0$ then $T^x(V_g)=V_{\;^{x}g}$ is also not zero and therefore $\opl_{h \in \co(g)}V_h$ is an equivariantized object of $\cc^F$. Since $V$ is a simple object it follows that $V$ is supported only on one orbit, namely the orbit $\co(g)$ of $g$. Thus $V=\oplus_{{h \in \co(g)}}V_{h}$. Moreover it follows that $V_g\in \cc_g^{F_g}$ and one can associate to the object $(V, {\mu^{x}_{V}}_{\{x \in F\}})\in \cc^F$ the pair $(g, M)$ with $M:=V_g$. Clearly $M$ is a simple object of $\cc_{g}^{F_{g}}$ if $V$ is a simple  object of $\cc^{F}$.
\md
Conversely, to any pair $(g, M)$ as above one associates the simple object $\Ind_{F_{g}}^{F}(M) \in \cc^{G}$.
\md
We have to show that the above two constructions are one inverse to the other.  It is easy to see that the component of order $g$ of $\Ind_{F_{g}}^{F}(M) \in \cc^{G}$ coincides to $M$ as an object of $\cc^{F_{g}}_{g}$.
\mdn Thus it remains to show that $V\simeq \ind^{F}_{F_{g}}(V_{g})$ if $V=\oplus_{g \in G}V_{g}$ is a simple object of $\cc^{F}$. 
This is equivalent with showing that
\bq
V\simeq \oplus_{r \in F/F_{g}}T^{r}(V_{g})
\eq
as objects of $\cc^{F}$.

Since $\mu_{V}^{x}:T^{x}(V) \xra{\simeq} V$ it follows that $\mu_{V}^{x}|_{T^{x}(V_{g})}:T^{x}(V_{g}) \xra{\simeq}V_{\;^{xg}}$ are also isomorphisms in $\cc$. This shows that $f: \ind_{F_{g}}^{F}(V_{g})\xra{\simeq} V$ is an isomorphism in $\cc^{F}$ where $f:=\oplus_{r \in F/F_{g}}\mu_{V}^{r}|_{T^{r}(V_{g})}$.
\epf
\br
Note that there is an embedding $\cc_{g}^{F_{g}}\subset \cc^{F_{g}}$ and one can apply the induced factor $\ind^{F}_{F_{g}}$ to $M$.
\er
\subsection{Definition of $S(g, M)$} Let $M \in \cc^{F_{g}}$ be a simple object.
Define by $S(g, M)$ as the simple induced object $\ind^{F}_{F_{g}}(M) \in \cc^F$ from above. Thus as objects of $\cc$ one has that $S(g, M):=\oplus_{r \in F/F_g}T^r(M)$ with the equivariant $F$-structure obtained from Equation \eqref{muind}.

\br\lb{shift} Note that if $(M , \mu_{M}) \in \cc_g^{F_g}$ then $(T^x(M), \;^x\mu_{M})\in \cc_{\;^xg}^{F_{\;^xg}}$ by Proposition \ref{conj}. Then item 2) of Proposition \ref{commir} implies that $S(g, M)\simeq S(\;^xg, T^x(M))$ for any $x \in F$, $g \in G$  and $M \in \cc_g^{F_g}$.
\er
\subsection{Tensor product formula}
Let $M\in \cc_g^{F_{g}}$ and $N \in \cc_h^{F_{h}}$ be two equivariant objects. Define the following equivariant object
\beq\label{mgh}
m_{g,h}(M, N):=\ind_{F_{g}\cap F_{h}}^{F_{gh}}(\res^{F_{g}}_{F_{g}\cap F_h}(M)\ot \res_{F_{g}\cap F_h}^{F_h}(N))\in \cc_{gh}^{F_{gh}}.
\eeq
\bt \lb{tpgr} Suppose that $F$ acts coherently on the $G$-graded fusion category $\cc$ with respect to a given action by group automorphisms of $F$ on $G$. With the above notations one has
\beq\lb{tepg}
S({g, M})\ot S({h,N})\simeq \opl_{x \in D}S(\;^xgh, \;\;m_{ _{\;^xg, h}}(T^x(M),N))
\eeq
where $D$ is a set of representative elements  for the double cosets $F_h \backslash F / F_g$.
\et
\bpf
One has by definition $S(g, M)=\ind_{F_g}^F(M)$. Applying formula (\ref{ti}) one has that
\beqarn
S({g, M})\ot S({h,N}) & = & \ind_{F_g}^F(M)\ot \ind_{F_h}^F(N)\\
& \simeq & \ind_{F_h}^F(\res_{F_h}^F(\ind_{F_g}^F(M)\ot N)).
\eeqarn

On the other hand applying Theorem \ref{macky} one has that 
\beqarn
\res_{F_h}^F(\ind_{F_g}^F(M))\simeq \oplus_{x \in D}\ind_{{F_h} \cap \;^x{F_g}}^{F_h}(\res^{\;^x{F_g}}_{\;^x{F_g}\cap {F_h}}(T^x(M)))
\eeqarn
Then applying Equation \eqref{ti2} one obtains that
\beqarn
\res_{F_h}^F(\ind_{F_g}^F(M))\ot N\simeq \oplus_{x \in D}(\ind_{{F_h} \cap \;^x{F_g}}^{F_h}(\res^{\;^x{F_g}}_{\;^x{F_g}\cap {F_h}}(T^x(M))\ot N) \\ \simeq \oplus_{x \in D}(\ind_{{F_h} \cap \;^x{F_g}}^{F_h}(\res^{\;^x{F_g}}_{\;^x{F_g}\cap {F_h}}(T^x(M)\ot \res^{{F_h}}_{\;^x{F_g}\cap {F_h}}(N))
\eeqarn
and therefore 
\beqn S({g, M})\ot S({h,N}) 
\simeq  \oplus_{x \in D}\ind_{\;^{x}F_{g}\cap F_{h}}^{F}(\res^{\;^{x}F_{g}}_{\;^{x}F_{g}\cap F_h}(T^{x}(M))\ot \res_{\;^{x}F_{g}\cap F_h}^{F_h}(N))
\eeqn
which by definition coincides to $\opl_{x \in D}S(\;^xgh, \;\;m_{ _{\;^xg, h}}(T^x(M),N))$.
\epf
\bc Suppose that $F$ acts coherently on the $G$-graded fusion category $\cc$ with respect to a given action by group automorphisms of $F$ on $G$. With the above notations it follows that
\beq\lb{ditf}
\cc^F \simeq \oplus_{g \in \Gm}\cc_g^{F_g}
\eeq
as indecomposable $\cc_{1}^F$-bimodule categories, where $\Gm$ is a set of representative elements for the orbits of the action of $F$ on $G$ .\ec
\bpf
Remark that $\cc_{1}^{F}$ is a tensor subcategory of $\cc^{F}$ consisting on those objects of $\cc^{F}$ supported only on $\cc_{1}$. Define the functor $F: \cc^{F} \ra \oplus_{g \in \Gm}\cc_g^{F_g}$ by sending $S(g, M)\mapsto M$. Note that formula (\ref{tepg}) implies that 
\beqn
S(1, M)\ot S(h, N)\simeq S(h, \res^{F}_{F_{h}}(M)\ot N)
\eeqn
and
\beqn
S(h, N)\ot S(1, M)\simeq S(h, N\ot \res^{F}_{F_{h}}(M))
\eeqn
which shows that each $\cc^{F_{h}}$ is a $\cc_{1}^{F}$-bimodule category. The above formulae also show that $\cc_{h}^{F_{h}}$ is an indecomposable $\cc_{1}^{F}$-bimodule category for all $h \in H$.
\epf
\subsection{On the Grothendieck ring of an equivariantization under a coherent action}\lb{withsp}In this subsection we show that the Grothendieck ring of an equivariantization under a coherent action has the structure of a Green ring as introduced in \cite{scoh}.
\subsubsection{On the rings introduced by Witherspoon and Bouc} In this subsection we recall the Green rings introduced in \cite{scoh}. 

Let $F$ be a finite group acting by group automorphisms on another finite group $G$. Suppose  that 
$A=\oplus_{g \in G}A(g)$ is a graded vector space endowed with two linear structures: $m_{g,h}:A(g)\ot A(h)\ra A(gh)$
and  $c_{x,g}:A(g)\ra A(\;^xg)$ satisfying the following compatibilities:

$(\mathrm{C1})$
$c_1=\id \text{and\;} c_{xy}=c_x c_y$ where 
 $c_x:=\oplus_{g \in G}\;c_{x,g}:A\ra A$.

$(\mathrm{C2})$ $c_{x, g}=\id_{A(g)}$ if $x \in F_{g}$.

$(\mathrm{C3})$ $c_xm_{g,h}=m_{\;^xg,\;^xh}(c_x\times c_x)
$

$(\mathrm{C}4)$ There is an element $1 \in A(1)$ such that
$
c_x(1)=1 \; \text{for all\;} x \in L \; \text{and\;} m_{1,g}(1, \al_g)=m_{g,1}(\al_g,1)=\al_g
$.

Let $A^F:=\{a \in A\;|\; c_x(a)=a\;\text{for all}\;x\in F\}$
be the subspace of $F$-invariants elements of $A$.

$(\mathrm{C5})$ For any $g \in G$ and $\al, \beta, \gamma \in A^{F}$ one has that
\beqn
\sum_{(d,e,f)\in T_g}m_{de, f}(m_{d,e}(\al_d, \beta_e),\gm_f)=\sum_{(d,e,f)\in T_g}m_{d, ef}(\al_d, m_{e,f}(\beta_e,\gm_f))
\eeqn
where the set $T_g$ is defined as follows. Note that the stabilizer subgroup $F_g$ acts diagonally on the set $\{(d,e,f)\in G\times G\times G\;|\; def=g\}$. By $(\mathrm{C}1)$-$(\mathrm{C}4)$ the left and right members of the previous equality do not depend on the chosen set $T_g$ of representative elements for the orbits of the action of $F_g$ on the above set.

Define a multiplicative structure on $A^{F}$ given for $\al, \beta \in A^{F}$ by 
\beq\tag{M}\label{multo}
(\al\beta)_g=\sum_{\{(h,k)\in F_g/G \times G\;| \;hk=g\}}m_{h,k}(\al_h,\beta_k)
\eeq
where the action of $F_g$ on $G\times G$ is diagonal.
Then it is shown in \cite{scoh} that under the conditions $(\mathrm{C}1)$-$(\mathrm{C}4)$ the above multiplication on $A^F$ is associative if and only if condition $(\mathrm{C}5)$ is also satisfied. Moreover, this multiplication not depend on the choice of representative set $T_{g}$.
\bt\lb{idr}
Let $F$ be a finite group acting coherently on a $G$-graded fusion category $\cc$ with respect to a given action by group automorphisms of $F$ on $G$. Then the Grothendieck ring of $\cc^G$ has the multiplicative structure of Equation \eqref{multo}.
\et
\bpf
Let $A=\oplus_{g \in G} A(g)$ where $A(g)=K_0(\cc^{F_g})$. Using Proposition \ref{conj} one can define $c_{x, g}: A(g)\ra A(\;^xg)$ by $[M]\mapsto [T^x(M)]$. Define also $m_{g,h}:A(g)\times A(h) \ra A(gh)$ via the map $m_{g,h}$ from Equation \eqref{mgh}. Then it is easy to verify that the compatibility conditions $(\mathrm{C}1)$,  $(\mathrm{C2})$ and $(\mathrm{C4})$ from the previous subsection are satisfied. Condition $(\mathrm{C3})$ is verified in the lemma below.
Moreover it is clear that $K_0(\cc^F)\hookrightarrow A^F$ via $[S(g, M)] \mapsto \oplus_{r \in F/F_{g}}[T^{r}(M)]$ defines an inclusion of vector spaces.

Using the description of simple modules of equivariantization from \cite{buna} it follows that in fact  $K_0(\cc^F)= A^F$. Indeed, by \cite[Remark 3.12]{buna} a $k$-linear basis of $K_{0}(\cc^{F})$ is given by the elements  $\sum_{r \in F/F_{Y}}T^{r}(Y)$ where $Y$ runs through all the orbits of the isomorphisms classes of simple objects of $\cc$. Recall that  $F_{Y}:=\{x \in F\;|T^{x}( Y)\cong Y\}$ is the inertia subgroup.

On the other hand if $M \in \cc_{g}^{F_{g}}$ is a simple object then by \cite[Theorem 2.12]{buna} it follows that $M\simeq \Ind^{F_{g}}_{F_{Y}}(Y \ot \pi)$, for a simple object $Y\in \cc$,  a constituent of $M$ and some projective representation $\pi$ of $(F_{g})_{Y}:=F_{Y}\cap F_{g}$. Since a $k$-linear bases of  $A^{F}$ is given by $\sum_{r\in F/F_{g}}T^{r}(M)$ where $M \in \cc_{g}^{F_{g}}$ is a simple object it follows from the above description of $M$ that the same vectors $\sum_{r \in F/F_{Y}}T^{r}(Y)$ form also a $k$-linear basis of $A^{F}$.

It remains to show that the multiplication from Theorem \ref{tpgr}  coincides to the multiplication described in Equation \eqref{multo}.  Once this is proven, it follows that condition $(\mathrm{C}5)$ is also satisfied since the multiplication in $K_{0}(\cc^{F})$ is associative.

Note that as in \cite{scoh} one has that $\co(g)\co(h)=\sqcup_{x \in F_{h }\slash F \backslash F_{g}}\co(\;^{x}gh)$. On the other hand multiplication formula from Theorem \ref{tpgr} shows
\beqarn
[S(g, M)][S(h, N)] & = & \sum_{x \in D}\sum_{r \in F/F_{\;^{x}gh}}[T^{r}(m_{^{x}g, h}(T^{x}(M), N))]\\ &=&  \sum_{x \in D}\sum_{r \in F/F_{\;^{x}gh}}m_{^{rx}g, \;^{r}h}([T^{rx}(M)], [T^{r}(N)])
\eeqarn
We have to show that this multiplication coincides to the one given in Equation \eqref{multo}.
For a fixed $x\in D$ note that $m_{^{rx}g, \;^{r}h}([T^{rx}(M)], [T^{r}(N)])$ with $r \in F/F_{\;^{x}gh}$ runs through all the orbit of of $\;^{x}gh$. Thus the term in $A(\;^{x}gh)$ of the above product coincides to $$\sum _{\{y \in D\:|\co(\;^{y}gh)=\co(\;^{x}gh)\}}m_{^{r_{y}y}g,\;^{r_{y}}h}([T^{r_{y}y}(M), T^{r_{y}}(N))$$ where $r_{y}\in F/F_{\;^{y}gh}$ is uniquely chosen such that $\;^{r_{y}y}g\;^{r_{y}}h=\;^{x}gh$.

One needs to show that the set $\{y \in D\:|\co(\;^{y}gh)=\co(\;^{x}gh)\}$ has the same cardinality as the set $[(\co(g)\times \co(h))\cap \{(a, b)\in G\times G\:|ab=\;^{x}gh\}]/F_{\;^{x}gh}$. In order to do this we construct a bijection between these two sets. If $\;^{r}g\;^{s}h=\;^{x}gh$ then send the orbit of the pair $(\;^{r}g,\;^{s}h)$ to the double coset $F_{h}s^{-1}rF_{g}$. Clearly this map is well defined. Conversely, define $F_{h}yF_{g}\mapsto \co((\:^{ry}g, \;^{r}h))$ where $r$ is chosen such that  $\:^{r}(\;^{y}gh)=\;^{x}gh$. It is easy to check that these two maps are one inverse to another.
\epf

\bl Suppose that $F$ acts coherently on a $G$-graded fusion category $\cc$ with respect to a given action by group automorphisms of $F$ on $G$.
Then with the above notations it follows that
\beq
T^{x}(m_{g, h}(M, N))\simeq m_{\;^{x}g, \;^{x}h}(T^{x}(M), T^{x}(N))
\eeq
for all $M\in \cc^{F_{g}}$ and $N\in \cc^{F_{h}}$.
\el
\bpf
As objects of $\cc$ one has that \beqn m_{g, h}(M, N)=\oplus_{r \in F_{gh}/F_{g}\cap F_{h}}T^{r}(M\ot N).\eeqn On the other hand \beqn m_{\;^{x}g, \;^{x}h}(T^{x}(M), T^{x}(N))=\oplus_{r \in F_{gh}/F_{g}\cap F_{h}}T^{xr^{-1}x}(T^{x}(M)\ot T^{x}(N)).\eeqn as objects of $\cc$.

It can be checked directly that $F: T^{x}(m_{g, h}(M, N))\ra m_{\;^{x}g, \;^{x}h}(T^{x}(M), T^{x}(N))$ given on components by
\beqn
T^{x}(T^{r}(M\ot N)) \xra{(T_{2}^{x,r})_{M\ot N}} T^{xr}(M\ot N)\xra{(T_{2}^{xrx^{-1},x})^{-1}_{M\ot N}} \eeqn \vskip -0,5cm \beqn \xra{(T_{2}^{xrx^{-1},x})^{-1}_{M\ot N}} T^{xrx^{-1}}(T^{x}(M\ot N))\xra{T^{xrx^{-1}}((T_{2}^{x})^{M, N})} T^{xrx^{-1}}(T^{x}(M)\ot T^{x}(N))
\eeqn
is an isomorphism in $\cc^{F_{\;^{x}(gh)}}$.

If $z' \in F_{gh}$ note that the equivariant structure of $m_{g, h}(M, N)$ is given on the components by 
\beqn
\mu^{z, r}_{m_{g, h}(M, N)}:T^{z}(T^{r}(M\ot N))\xra{(T_{2}^{z, r})_{M\ot N}} T^{zr}(M\ot N) \xra{(T_{2}^{r', h})^{-1}_{M\ot N}}\eeqn \vskip -0,5cm \beqn  T^{r'}(T^{l}(M\ot N))
\xra{T
^{r'}((T^{l}_{2})^{M, N})} T^{r'}(T^{l}(M)\ot T^{l}(N))\xra{T^{r'}(\mu_{M}^{l}\ot \mu_{N}^{l})} T^{r'}(M\ot N) \eeqn \vskip -0,5cm \beqn 
\eeqn
where $zr=r'l$ with $l \in F_{g}\cap F_{h}$ and $r' \in F_{gh}/F_{g}\cap F_{h}$.
Using Equation \eqref{conjmu} it follows that the equivariant structure of $T^{x}(m_{g, h}(M, N))$ is given on components by
\beqn
\;^{x}\mu^{xzx^{-1},\; r}_{T^{x}(m_{g, h}(M, N))}:T^{xzx^{-1}}(T^{x}(T^{r}(M\ot N))\xra{(T_{2}^{xzx^{-1}, x})_{M\ot N}} T^{xz}(T^{r}(M\ot N))\xra{} \eeqn \vskip -0,5cm \beqn \xra{(T_{2}^{x, z})^{-1}_{M\ot N}} T^{x}(T^{z}(T^{r}(M\ot N)) \xra{T^{x}(\mu^{z, r'}_{m_{g, h}(M, N)})}T^{x}(T^{r'}(M\ot N)).
\eeqn
On the other hand the equivariant structure $\nu^{xzx^{-1},\; r}_{ _{m_{\;^{x}g, \;^{x}h}(T^{x}(M), T^{x}(N))}}$ of the object $m_{\;^{x}g, \;^{x}h}(T^{x}(M), T^{x}(N))$
is given on the components by
\beqn
T^{xzx^{-1}}(T^{xr^{-1}x}(T^{x}(M)\ot T^{x}(N))) \xra{(T_{2}^{xzx^{-1}, xrx^{-1}})_{M\ot N}} T^{xzrx^{-1}}(T^{x}(M)\ot T^{x}(N))
\xra{} 
\eeqn
\beqn \xra{(T_{2}^{x, z})^{-1}_{M\ot N}} T^{x}(T^{z}(T^{r}(M\ot N)) 
\xra{T^{x}(\mu^{z, r'}_{m_{g, h}(M, N)})}T^{x}(T^{r'}(M\ot N)).
\eeqn
Therefore verifying that $F$ is an isomorphism resumes to the commutativity of the following diagram (D1) made of solid arrows below. The commutativity of the bottom left rectangle of diagram (D1) follows from commutativity of diagram (D2). For shortness the maps in the diagrams are omitted but they are all uniquely determined from the group action of $G$ on $\cc$.
\subsubsection{Grothendieck rings of abelian cocentral extensions} In \cite[Theorem 4.8]{scoh} it is shown that the Grothendieck rings $\mathcal{G}_{0}(H)$ associated to abelian cocentral extensions  have the multiplication structure given in Equation \eqref{multo}. Proposition \ref{cocsp} implies that the same result holds for any cocentral extension of semisimple Hopf algebras.

\br Note that compatibility condition $(\mathrm{C2})$ is not stated in \cite{scoh} on page 5 although it is stated as a property for Grothendieck groups of cocentral extensions on the last page of the paper. \er

\subsubsection{On the Grothendieck ring of the center of a fusion category} Suppose that $\cc$ is a $G$-graded fusion category $\cc=\oplus_{g \in G}\;\cc_{g}.$
Then by \cite[Theorem 4.1]{gnn} its Drinfeld center $\cz(\cc)\simeq \cz_{\cc_{1}}(\cc)^{G}$, the equivariantization of the relative center $\cz_{\cc_{1}}(\cc)$ by a certain action of the finite group $G$. Moreover \cite[Theorem 3.2]{gnn} shows that the relative center $\cz_{\cc_{1}}(\cc)$ is a $G$-crossed braided fusion category. In view of Example \ref{gbb} one can apply Theorem \ref{idr}. It follows that the Grothendieck ring of $\cz(\cc)$ has the ring structure described in Equation \eqref{multo}. \bll{Note that for the Grothendieck ring of a Drinfeld double of a semisimple Hopf algebra a similar description was already obtained in \cite{gdd}. }
\begin{center}
\begin{equation}\tag{Diagram D1}
{\Tiny
\begin{tikzpicture}
  \matrix (m) [matrix of math nodes, row sep=5. 5 em,column sep=1em,minimum width=2em]
  {
     T^{xzx\inv}(T^{x}(T^{r}(M\ot N))) & 
     T^{xzx\inv}(T^{xr}(M\ot N)) &   
     T^{xzx\inv}(T^{xrx^{-1}}(T^{x}(M\ot N))) \\
    T^{xz}(T^{r}(M\ot N)) & & T^{xzx\inv}(T^{xrx^{-1}}(T^{x}(M)\ot T^{x}(N))) \\
    T^{xzr} (M\ot N) & T^{xzrx^{-1}} (T^{x}(M\ot N))&  T^{xzrx^{-1}} (T^{x}(M)\ot T^{x}(N)) \\
     T^{x}(T^{zr} (M\ot N) )&  T^{xr'x^{-1}} (T^{xlx^{-1}}(T^{x}(M\ot N))  & T^{xr'x^{-1}} (T^{xlx^{-1}}(T^{x}(M)\ot T^{x}(N))) \\
     T^{x}(T^{r'}(T^{l}(M\ot N)))&  T^{xr'x^{-1}}(T^{xl}(M\ot N))  &   T^{xr'x^{-1}}(T^{xlx^{-1}}(T^{x}(M))\ot T^{xlx^{-1}}(T^{x}(N)))\\
    T^{x}(T^{r'}(T^{l}(M)\ot T^{l}(N))) &  T^{xr'x^{-1}}(T^{x}(T^{l}(M\ot N))) & T^{xr'x^{-1}}(T^{xl}(M)\ot T^{xl} (N))   \\
     T^{x}(T^{r'}(M\ot N))) &   T^{xr'x^{-1}}(T^{x}(T^{l}(M)& T^{xr'x^{-1}}(T^{x}(T^{l}(M)\ot T^{x}(T^{l} (N))))\\
     T^{xr'}(M\ot N)& T^{xr'x^{-1}}(T^{x}(M\ot N)) & T^{xr'x^{-1}}(T^{x}(M)\ot T^{x}(N)) .\\
     };
     \path[-stealth]
     (m-1-1) edge node [above] {$$}(m-1-2)
     (m-1-2) edge node [above] {$$}(m-1-3) 
     (m-8-1) edge node [below] {$$}(m-8-2)
      (m-8-2) edge node [below] {$$}(m-8-3)
        (m-3-1) edge [dashed] node  [above] {$$}(m-3-2)
        (m-3-2) edge [dashed] node [above] {$$}(m-3-3)
        (m-4-2) edge [dashed] node [above] {}(m-4-3)
        (m-5-2) edge [dashed] node [above] {}(m-6-3)
        (m-7-2) edge [dashed] node [above] {}(m-7-3)
      (m-8-2) edge  node [above] {}(m-8-3)
       (m-1-1) edge node [left] {$$} (m-2-1)
       (m-2-1) edge node [left] {$$} (m-3-1)
       (m-3-1) edge node [left] {$$} (m-4-1)
       (m-4-1) edge node [left] {$$} (m-5-1)
       (m-5-1) edge node [left] {$$} (m-6-1)
     (m-6-1) edge node [left] {$$} (m-7-1)
     (m-7-1) edge node [left] {$$} (m-8-1)
     (m-1-3) edge node [right] {$$} (m-2-3)
       (m-2-3) edge node [right] {$$} (m-3-3)
       (m-3-3) edge node [right] {} (m-4-3)%
       (m-4-3) edge node [left] {$$} (m-5-3)%
       (m-5-3) edge node [left] {} (m-6-3)
     (m-6-3) edge node [left] {$$} (m-7-3) 
     (m-7-3) edge node [left] {$$} (m-8-3) 
     (m-1-2) edge [dashed] node [right] {$$} (m-3-1)
        (m-1-3) edge [dashed] node [left] {$$} (m-3-2)
        (m-3-2) edge [dashed] node [left] {$$} (m-4-2)
          (m-4-2) edge [dashed] node [left] {$$} (m-5-2)
          (m-5-2) edge [dashed] node [left] {$$} (m-6-2)
           (m-6-2) edge [dashed] node [left] {$$} (m-7-2)
            (m-7-2) edge [dashed] node [left] {$$} (m-8-2)
            ;
  \end{tikzpicture}
}
\end{equation}
\end{center}

{\Tiny
\begin{center}
\begin{equation}\tag{Diagram D2}
\begin{tikzpicture}
  \matrix (m) [matrix of math nodes, row sep=4.5 em,column sep=1 em, minimum width=2em]
  {
    T^{xzr} (M\ot N) & & T^{xzrx^{-1}} (T^{x}(M\ot N)) \\
     T^{x}(T^{zr} (M\ot N) )& & T^{xr'x^{-1}} (T^{xlx^{-1}}(T^{x}(M\ot N))  \\
     T^{x}(T^{r'}(T^{l}(M\ot N)))&  T^{xr'}(T^{l}(M\ot N)) & T^{xr'x^{-1}}(T^{xl}(M\ot N)) \\
    T^{x}(T^{r'}(T^{l}(M)\ot T^{l}(N))) & T^{xr'}(T^{l}(M)\ot T^{l}(N) )& T^{xr'x^{-1}}(T^{x}(T^{l}(M\ot N)))   \\
     T^{x}(T^{r'}(M\ot N))) & &  T^{xr'x^{-1}}(T^{x}(T^{l}(M)\\
     T^{xr'}(M\ot N)& & T^{xr'x^{-1}}(T^{x}(M\ot N)) \\
     };
      \path[-stealth]
     (m-1-1) edge node [left] {$$} (m-2-1)
       (m-2-1) edge node [left] {$$} (m-3-1)
       (m-3-1) edge node [left] {$$} (m-4-1)%
       (m-4-1) edge node [left] {$$} (m-5-1)%
       (m-5-1) edge node [left] {$$} (m-6-1)%
        (m-1-3) edge node [left] {$$} (m-2-3)
       (m-2-3) edge node [left] {$$} (m-3-3)
       (m-3-3) edge node [left] {$$} (m-4-3)%
       (m-4-3) edge node [left] {$$} (m-5-3)%
       (m-5-3) edge node [left] {$$} (m-6-3)%
        (m-1-1) edge node [above] {$$}(m-1-3)
        (m-6-1)edge node [above] {$$}(m-6-3)
         (m-1-1) edge [dashed] node [above] {$$}(m-3-3)
         (m-1-1) edge [dashed] node [above] {$$}(m-3-2)
         (m-3-1) edge [dashed] node [above] {$$}(m-3-2)
         (m-3-2) edge [dashed] node [above] {$$}(m-4-3)
         (m-4-1) edge [dashed] node [above] {$$}(m-4-2)
         (m-3-2) edge [dashed] node [above] {$$}(m-4-2) 
         (m-4-2) edge [dashed] node [above] {$$}(m-5-3)
         (m-4-2) edge [dashed] node [above] {$$}(m-6-1)           ;
  \end{tikzpicture}
  \end{equation}
  \end{center}
}
\epf

\bibliographystyle{amsplain}
\bibliography{gfr}
\newpage
\bne
\item \bll{What if passing ot the Lrauer group?}
 \mdn
   \bll{need to work with arbitrary representative elements  for cosets since maximal chains might intersect}
\bll{ The compatibility \eqref{P} follows form equation \eqref{ro-2}. Equations \eqref{idS}, \eqref{Sid}  and \eqref{unul} are clear.
}
\item    \mdn
The compatibilities between the natural transformations follow from the following general Lemma:\mdn
We prove relation \eqref{N}. \bl
   Suppose that $z=a_{1}$
   \el

\bll{DEFINITION OF $Q^{L}_{H,a}$ also depends on some chosen representative elements }.
\bll{ Note that in the above Mackey functor one has equality in $\mathrm{(CM1)}$.}

\item \bll{Abelian versus $k$-linear.}
\item work with two group actions instead.
\item ask Adriana how to use landscape
\item the first part works for tensor categories.
\item For induction one has that isomorphisms not equalities since the repr of cosets change. 
\item morphisms between $\cc^{G}$-module functors.
\item module structure of a $\cc^{G}$-module functor.
\ene
\end{document}